\documentclass[leqno,11pt]{amsart}

\usepackage{geometry}

\usepackage[all]{xy} 

\usepackage{pgfplots}

\usepackage{mathtools}
\usepackage{amsmath, amssymb, amsfonts, latexsym, mdwlist, amsthm}
\usepackage{subfig}
\usepackage{graphicx}
\usepackage{wrapfig}
\usepackage{comment}
\usepackage{tikz-cd}

\usepackage[bookmarks, colorlinks, breaklinks, pdftitle={},
pdfauthor={Soheyla Feyzbakhsh}]{hyperref}
\hypersetup{linkcolor=blue,citecolor=blue,filecolor=black,urlcolor=blue}

\usepackage{tikz}
\usetikzlibrary{calc,trees,positioning,arrows,chains,shapes.geometric,%
    decorations.pathreplacing,decorations.pathmorphing,shapes,%
    matrix,shapes.symbols}

\tikzset{
>=stealth',
  punktchain/.style={
    rectangle,
    rounded corners,
    draw=black, thick,
    minimum height=3em,
    text centered,
    on chain},
  line/.style={draw, thick, <-},
  element/.style={
    tape,
    top color=white,
    bottom color=blue!50!black!60!,
    minimum width=8em,
    draw=blue!40!black!90, very thick,
    text width=10em,
    minimum height=3.5em,
    text centered,
    on chain},
  every join/.style={->, thick,shorten >=1pt},
  decoration={brace},
  tuborg/.style={decorate},
  tubnode/.style={midway, right=2pt},
}

\usepackage{paralist}
\setdefaultenum{(a)}{(i)}{}{}
\usepackage{enumitem} 
\usepackage{graphicx}

\makeatletter
\newtheorem*{rep@theorem}{\rep@title}
\newcommand{\newreptheorem}[2]{%
\newenvironment{rep#1}[1]{%
 \def\rep@title{#2 \ref{##1}}%
 \begin{rep@theorem}}%
 {\end{rep@theorem}}}
\makeatother

\newtheorem{Thm}{Theorem}[section]
\newreptheorem{Thm}{Theorem}
\newtheorem{Prop}[Thm]{Proposition}

\newtheorem{Lem}[Thm]{Lemma}

\newtheorem{Cor}[Thm]{Corollary}
\newreptheorem{Cor}{Corollary}

\newreptheorem{Con}{Conjecture}

\newtheorem{thm-int}{Theorem}

\theoremstyle{definition}
\newtheorem{Def-s}[Thm]{Definition}

\newtheorem{Rem}[Thm]{Remark}

\newcommand{\ignore}[1]{}
\usepackage{amsmath,calligra,mathrsfs}

\newcommand{\Ale}[1]{\textcolor{blue}{#1}}

\renewcommand\_{^{}_}

\newcommand\C{\mathbb C}

\newcommand\Z{\mathbb Z}
\newcommand\cA{\mathcal A}
\newcommand\cD{\mathcal D}
\newcommand\cH{\mathcal H}
\newcommand\cT{\mathcal T}

\newcommand\rd{\operatorname{\textbf{r}}}
\newcommand\dd{\operatorname{\textbf{d}}}
\newcommand\nd{\operatorname{\textbf{n}}}

\newcommand\T{{\mathcal T}_C}
\newcommand\cO{{\mathcal O}_C}

\newcommand\Hom{\operatorname{Hom}}
\newcommand\Ext{\operatorname{Ext}}
\newcommand\ext{\operatorname{ext}}

\newcommand\Coh{\operatorname{Coh}}
\newcommand\cl{\operatorname{cl}}

\newcommand\Ho{\mathscr{H}\text{\kern -3pt {\calligra\large om}}}

\newcommand\rk{\operatorname{rk}}

\begin{document}
	\title[Butler's conjecture]
	{Butler's conjecture for K3 curves\vspace{-2mm}}

\author{Soheyla Feyzbakhsh}
\address{Department of Mathematics, Imperial College, London SW7 2AZ, United Kingdom}
\email{s.feyzbakhsh@imperial.ac.uk}
    
\author{Aliaksandra Novik}
\address{Department of Mathematics, Imperial College, London SW7 2AZ, United Kingdom}
\email{a.novik24@imperial.ac.uk}

\author{Angela Ortega}
\address{Institut für Mathematik, Humboldt Universität zu Berlin,  Unter den Linden 6, D-10099 Berlin, Germany    }
\email{ortega@math.hu-berlin.de}

\begin{abstract} \noindent
Butler's conjectures predict that, on a general smooth curve $C$ of genus $g\ge3$, a general coherent system $[\mathcal{O}_C\otimes V\xrightarrow{\varphi}E]$ with semistable underlying bundle $E$ and $\dim V>\operatorname{rk}(E)$ has surjective evaluation map $\varphi$ and stable kernel.

In this paper, we prove Butler’s conjectures for $g\ge8$ in the two cases in which the Brill--Noether geometry of a general curve is known, namely rank one and the range $\mu(E)>2g-2$. By deformation and openness, it suffices to establish the result on a single Brill--Noether--Petri general curve. We do so for curves on K3 surfaces of Picard rank one, using the wall-crossing theory of coherent systems developed by the first two authors.

\end{abstract}
	\maketitle	



\section{Introduction}
Let $C$ be a smooth complex projective curve. For positive integers $r$, $d$ and $n$, we consider coherent systems
\[
  T=\bigl[\mathcal{O}_C\otimes V\xrightarrow{\ \varphi\ }E\bigr]
\]
such that:
\begin{enumerate}
  \item $V$ is a vector space of dimension $n$;
 \item the induced map $H^0(\varphi)\colon V\to H^0(C,E)$ is injective;
  \item $E$ is a semistable vector bundle of rank $r$ and degree $d$; such that for every nonzero proper subbundle \(E'\subset E\) with \(\mu(E')=\mu(E)\) we have 
  \[
\frac{\dim\bigl(V\cap H^0(C,E')\bigr)}{\operatorname{rk}E'}
<
\frac{\dim V}{\operatorname{rk}E}
\]
\end{enumerate}
We denote their moduli space by $M_\infty(r,d,n)$. 
Butler's conjecture, proposed in 1997, predicts the following: if $C$ is a general curve of genus $g\ge 3$, $n>r$, and $M_\infty(r,d,n)$ is nonempty, then a general system
\[
  T=\bigl[\mathcal{O}_C\otimes V\xrightarrow{\ \varphi\ }E\bigr]\in M_\infty(r,d,n)
\]
has surjective evaluation map $\varphi$ and slope-stable kernel $\ker\varphi$, see \cite[Conjectures~1 and~2]{Butler1997}.

Butler's conjecture has been studied extensively since its formulation; see Section~\ref{work} for related work. Its geometric appeal lies in the dual span construction, which predicts birational correspondences between moduli spaces of coherent systems of different ranks, thereby relating distinct Brill--Noether loci on the same curve.

In this paper, we address the conjecture in rank one, and in arbitrary rank with slope greater than $2g-2$: two ranges in which the relevant Brill--Noether theory is known. We focus on smooth curves on general K3 surfaces. Such curves satisfy the classical Brill--Noether theorem: their loci of linear series have the same nonemptiness and dimension properties as those of a general curve \cite{lazarsfeld:brill-noether}. Our main theorem is the following. 

\begin{Thm}\label{thm-main}
Let \(X\) be a smooth K3 surface with \(\operatorname{Pic}(X)=\mathbb{Z}H\), and let \(C\in|H|\) be a general smooth curve of genus \(g\ge 8\). Suppose that one of the following holds:
  \begin{itemize}
    \item $r=1$ and $n\ge 2$; or
    \item $n>r>0$ and $d>r(2g-2)$.
  \end{itemize}
  If $M_\infty(r,d,n)$ is nonempty, then a general coherent system
  \[
    T=\bigl[\cO\otimes V\xrightarrow{\ \varphi\ }E\bigr]\in M_\infty(r,d,n)
  \]
  has surjective evaluation map $\varphi$, and $\ker\varphi$ is slope-stable.
\end{Thm}

For each fixed numerical type $(r,d,n)$ covered by Theorem \ref{thm-main}, the same conclusion holds on a general smooth curve $C$ of genus $g\ge 8$. 
In the high-slope case, $E$ is slope-semistable by the definition of $M_\infty(r,d,n)$. Since $\operatorname{Ext}^2(E,E)=0$ on a curve, the bundle deforms with the curve, and since $H^1(C,E)=0$, its chosen sections lift. In rank one, we may start with a general smooth member of $|H|$, where Petri injectivity \cite{lazarsfeld:brill-noether} gives the same deformation property for linear series. Surjectivity of the evaluation map and slope stability of its kernel are open conditions. Since the relevant relative moduli space is smooth over the base at the points constructed above, the locus where both conditions hold dominates the moduli of curves, and hence is nonempty over a general curve. It is dense by irreducibility in the high-slope case, and in rank one when $\rho>0$. When $\rho=0$, every point of the finite fibre over the chosen K3 curve satisfies the conclusion. Properness of the relative space of linear series then gives the same conclusion for every linear series on a nearby curve.

\medskip

Theorem \ref{thm-main} follows from Theorems \ref{thm-any rank-high slope} and \ref{thm-rank one}, which treat the high-slope and rank-one cases, respectively. In fact, the former establishes a slightly more general version of the high-slope statement.

As a direct consequence of the rank one case, we obtain the following corollary. 

\begin{Cor}\label{cor:bundles}
Under the assumptions on $C$ in Theorem~\ref{thm-main}, a
globally generated stable bundle of rank $R\ge1$ and degree
$d>0$ exists if and only if
\begin{equation}\label{eq:corollary}
\begin{cases}
\rho(g,R,d):=g-(R+1)(g-d+R)\ge0,& d\le2g-2,\\
R\le d-g,& d>2g-2.
\end{cases}
\end{equation}
\end{Cor}

\medskip

\subsection{Strategy of the proof} 

We work with the two-dimensional slice of Bridgeland stability conditions $\sigma_{b,w}$ on the bounded derived category of coherent systems constructed in \cite{FN}. The large volume chamber recovers slope stability of the underlying bundle $E$. For a fixed class $v=(r,d,n)$, the stable object
\[
  \bigl[\mathcal{O}_C\xrightarrow{\ \mathrm{id}\ }\mathcal{O}_C\bigr]
\]
determines the Brill--Noether wall.
If a coherent system $T=[\mathcal{O}_C\otimes V\xrightarrow{\ \varphi\ }E]$ remains stable down to this wall, then $\varphi$ is surjective and its kernel is slope-semistable. The problem therefore reduces to showing that a general system reaches this wall without becoming unstable.

After applying the degree-shift Proposition \ref{prop:degree-shift-gen} to reduce the numerical range, we analyse the possible first destabilizing walls. For each wall, we estimate the dimension of the destabilized locus by deforming the two factors and their extension. The key numerical estimates control the obstruction spaces through bounds on sections and on the Harder--Narasimhan filtrations of the underlying bundles. This treats both sheaf-theoretic and non-sheaf-theoretic walls, and shows that each destabilized locus has dimension strictly smaller than that of the large volume moduli space.

Local finiteness of walls along the chosen compact path then allows us to avoid all these loci simultaneously. Consequently, a general coherent system remains semistable down to the Brill--Noether wall, which yields global generation and semistability of its evaluation kernel.

\subsection{Related work} \label{work}
Butler himself proved that the complete evaluation kernel $M_E = \ker (\cO \otimes H^0(E) \xrightarrow{\text{ev}} E)$ is semistable when $E$ is semistable and $\mu(E)\ge 2g$, and stable when $E$ is stable and $\mu(E)>2g$ \cite{Butler1994}. This extended earlier results for line bundles: Paranjape and Ramanan showed that the kernel bundle of the canonical series is semistable, and stable for non-hyperelliptic curves \cite{ParanjapeRamanan1988}, while Ein and Lazarsfeld proved that $M_L$ is stable when $\deg L>2g$ and semistable when $\deg L=2g$ \cite{EinLazarsfeld1992}. Butler's subsequent conjecture extended this question to general coherent systems and proposed birational correspondences induced by the dual span construction \cite{Butler1997}. Further results relating kernel stability to the Clifford index, linear stability and the codimension of the chosen space of sections were obtained by Camere, Mistretta, and Mistretta--Stoppino \cite{Camere2008,Mistretta2008,MistrettaStoppino2012}. Brambila-Paz and Ortega studied the corresponding Brill--Noether bundles on special curves \cite{BPO2007}.

The rank-one semistability conjecture was proved by Bhosle, Brambila-Paz and Newstead using wall-crossing for coherent systems: on a general complex curve, a general generated linear series of any admissible numerical type has semistable evaluation kernel \cite{BBPN2015}. More recently, Farkas and Larson proved stability of the restricted tangent bundle, and hence of the evaluation kernel, for general Brill--Noether curves in projective space, with a specific genus-two exception, in arbitrary characteristic \cite{FarkasLarson2025}; see also the survey of Coskun, Larson and Vogt \cite{CLV2026}.

In higher rank, the conjecture remains open in general. Bhosle, Brambila-Paz and Newstead studied coherent systems of type $(r,d,r+1)$ on Petri curves via the dual span construction \cite{BBPN2008}. Brambila-Paz, Mata-Guti\'errez, Newstead and Ortega established several cases, particularly for rank-two systems with four sections \cite{BPMGNO2019}. Castorena, Hitching and Luna extended Mumford's notion of linear stability to generated coherent systems of higher rank, and proved that linearly stable systems of type $(2,d,4)$ with $d$ sufficiently small have stable dual span bundle \cite{CHL2025}. Subsequent work of Castorena--Hitching and Bajravani--Ortega obtained results for rank-two systems with five sections \cite{CastorenaHitching2024,BajravaniOrtega2026}. Castorena and Hitching further showed that linear stability of $(E,V)$ is necessary, but not sufficient, for stability of the kernel bundle $M_{E,V}$ in higher rank \cite{CastorenaHitching2025}.

\subsection*{Acknowledgments}
We thank Arend Bayer, Gavril Farkas, Chunyi Li, Zhiyu Liu, Richard Thomas and Yukinobu Toda for very helpful discussions. S.F.\ acknowledges support from the Royal Society (URF/R1/23119).

\subsection*{Use of AI tools}
We used ChatGPT, Astra~6 and Claude Opus~5.5 to assist with the numerical analysis of wall-crossing in Sections~3 and~4. All results were checked by the authors, who take full responsibility for the content of this paper.

\bigskip

\section{Two-dimensional slice of Bridgeland stability conditions}\label{section.2}
In this section we give a brief review on the construction of the derived category of coherent systems on a smooth projective curve $C$ of genus $g\geq 2$ and its two-dimensional slice of Bridgeland stability conditions introduced in~\cite{FN}. We then recall the corresponding wall-and-chamber structure and the moduli spaces associated to the large-volume limit. 

A \textit{(generalised) coherent system} on $C$ is a triple $(V,E,\varphi)$, where $V$ is a finite-dimensional $\mathbb{C}$–vector space, $E$ is a coherent sheaf on $C$, and $\varphi\colon \cO\otimes V\to E$ is an arbitrary sheaf morphism.\footnote{In much of the literature on coherent systems one assumes $H^0(\varphi)$ is injective; in this paper we allow $\varphi$ to be arbitrary.}  We write such a triple as
\[
[\cO \otimes V \xrightarrow{\ \varphi\ } E].
\]
The category of these triples, denoted $\mathcal{T}_C$, is abelian, and we write $\mathcal{D}(\mathcal{T}_C)$ for its bounded derived category. 

The numerical Grothendieck group of $\cD(\cT_C)$ is identified with $\mathbb{Z}^3$ via the class map, which associates to any object $T=[\cO\otimes V\xrightarrow{\,\varphi\,} E]$ the vector
\[
\cl(T)
= \big(\rk E,\ \deg E,\ \dim_{\C} V\big)
\eqqcolon (\rd(T),\dd(T),\nd(T)).
\]
For two objects $T_i = [\cO \otimes V_i \xrightarrow{\varphi _i} E_i] \in \mathcal{D}(\mathcal{T}_C)$ for $i=1, 2$ with $\cl(T_i) = (r_i, d_i, n_i)$ the Euler pairing is given as
\begin{align*}
\chi(T_1, T_2) \coloneqq & \sum_{k}  (-1)^k\ext^k(T_1 , T_2)\\
= & r_1 (d_2 +r_2(1-g)) -d_1r_2 + n_1(n_2 -d_2-r_2(1-g)). 
\end{align*}

There is a natural semiorthogonal decomposition
	 \begin{align*}
	    \mathcal{D}(\mathcal{T}_C) =\ &  \langle i_*\cD(\mathcal{V}) ,\ j_*\mathcal{D}(C) \rangle, 
	 \end{align*}
     with $i_*V = [\cO \otimes V \to 0]$ and $j_*F=[0\to F]$. Moreover, the adjoints $i^* \dashv i_*$ and $j_* \dashv j^!$ for any $T=[\cO \otimes V \to E]\in \mathcal{D}(\T)$ are defined as $i^*T = V$, $j^!T = E$. 

The two-dimensional slice of the space of stability conditions on $\cD(\T)$ is obtained by tilting the natural heart $\T$ with respect to a torsion pair. Therefore, we first extend the classical notion of $\mu$-stability of sheaves on a curve to triples. Fix $\alpha \in \mathbb{R}_{\geq 0}$. For any object $T = [\cO \otimes V \to E] \in \T$, we define the slope 
\begin{equation*}
\mu_{\alpha}(T)\ \coloneqq\ \left\{\!\!\!\begin{array}{cc} \frac{\deg(E)}{\rk(E)} + \alpha \frac{\dim V}{\rk(E)}
& \text{if }\rk(E)\ne0, \\
+\infty & \text{if }\rk(E)=0. \end{array}\right.
\end{equation*}
We say that $T \in \mathcal{T}_C$ is \textit{$\mu_{\alpha}$-(semi)stable} if for all non-trivial subobject $0 \neq T' \subset T$ in $\T$, we have $\mu_{\alpha}(T') < (\leq)\ \mu_{\alpha}(T/T')$.

By~\cite[Proposition 2.3]{He1996}, every object $T\in\cT_C$ admits a unique Harder–Narasimhan filtration with respect to $\mu_{\alpha}$-stability. We focus on the case $\alpha=0$ and write $\mu:=\mu_{0}$ for simplicity. We denote by $\mu^{+}(T)$ and $\mu^{-}(T)$ for the maximal and minimal slopes occurring in the HN filtration of $T$, respectively. For $b\in \mathbb{R}$, let $\mathbb{T}^b$ and $\mathbb{F}^b$ be full subcategories in $\T$, defined as 
\[
\mathbb{T}^b=\{T\in \T: \mu^{-}(T)>b\}, \qquad \mathbb{F}^b=\{T\in \T: \mu^{+}(T)\leq b\}. 
\]
The pair $(\mathbb{T}^b, \mathbb{F}^b)$ is a torsion pair. Tilting $\T$ with respect to it gives the heart of a bounded $t$-structure:
\begin{equation*}
        \cA(b) \coloneqq \langle \mathbb{T}^b,\  \mathbb{F}^b[1] \rangle.
    \end{equation*}

The two-dimensional slice of stability conditions obtained by this construction, is controlled by the \textit{Brill-Noether function} $\Phi_C \colon \mathbb{R} \to \mathbb{R}$, which is defined as 
\begin{equation*}
    \Phi_C(x) \coloneqq \limsup_{\mu \to x} 
    \left\{ \frac{h^0(C, F)}{\rk(F)} \colon\ \text{$F \in \Coh(C)$ is semistable with slope $\mu(F) =\mu$}     \right\}. 
\end{equation*}
This is the smallest upper semi-continuous function satisfying 
\begin{equation*}
    \frac{h^0(F)}{\rk(F)} \leq \Phi_C(\mu(F))
\end{equation*}
for every semistable sheaf $F$ on $C$. Define an open subset
\[
U_C : = \{(b, w): w> \Phi_C(b)\} \subset \mathbb{R}^2.
\]

\begin{Thm}~\cite[Theorem 3.3]{FN}
    There is a two-dimensional continuous family of stability conditions parametrized by  $(b,w) \in U_C$ given by $(b,w) \mapsto \sigma_{b,w} \coloneqq (\cA(b), Z_{b,w})$ for the the group homomorphism 
 \begin{equation*}
  Z_{b,w} \colon \mathcal{N}(\cD(T)) \rightarrow \mathbb{C} \;\;\;,\;\;\; Z_{b,w}(T) = -\mathbf n(T)+w\mathbf r(T)+i (\mathbf d(T)-b\mathbf r(T)).
  \end{equation*}
\end{Thm}

For any non-zero object $T \in \cA(b)$, the associate slope function is 
  \begin{equation*}
      \nu_{b,w}(T) = \frac{\nd(T)-w \rd(T)}{\dd(T)-b  \rd(T)}, 
  \end{equation*}
where we write $\nu_{b,w}(T)=+\infty$ when the denominator vanishes. We say $0\neq T \in \cD(\cT_C)$ is \textit{$\sigma_{b,w}$-(semi)stable} if and only if
 $T[k]\in\cA(b)$ for some $k\in\Z$, and $\nu_{b,w}(T')< (\le)\nu_{b,w}\big(T[k]/T')$ for all non-trivial subobjects $T'\hookrightarrow T[k]$ in $\cA(b)$.

We next recall the wall-and-chamber structure of $U_C$. We plot the $(b,w)$-plane simultaneously with the image of the projection map
\begin{equation*}
\Pi \colon \mathcal{N}(\mathcal{D}(\T)) \rightarrow \mathbb{R}^2 \;\;\;,\;\;\; \Pi(r,d,n) = \left(\dfrac{d}{r} , \dfrac{n}{r} \right),  
\end{equation*}
defined for the classes of nonzero rank. For any class $v=(r,d,n)$ of nonzero rank, the tilt slope is precisely the slope of the line joining the point $(b,w)$ to its projection $\Pi(v)$, with the convention that vertical lines have infinite slope.

\begin{Prop}~\cite[Proposition 6.1]{FN}
	Fix $v = (r, d, n)\in \mathcal{N}(\mathcal{D}(\T))$. 
    There exists a set of line segments $\{\ell_i\}_{i \in I}$ in $U_{C}$ (called ``\emph{walls}") which are locally finite and satisfy 
	\begin{itemize*}
	    \item[\emph{(}a\emph{)}] If \(r\neq 0\), then the line containing \(\ell_i\) passes through \(\Pi(v)\).
	    \item[\emph{(}b\emph{)}] If $r=0$ then all $\ell_i$ are parallel of slope $\frac{n}{d}$.
        \item [\emph{(}c\emph{)}] The line segments \(\ell_i\) terminate on the boundary \(\partial U_C\).
	   		\item[\emph{(}d\emph{)}] The $\sigma\_{b,w}$-(semi)stability of any $T \in \cD(\cT_C)$ of class $v$ is unchanged as $(b,w)$ varies within any connected component (called a ``\emph{chamber}") of $U_C \setminus \bigcup_{i \in I}\ell_i$.
		\item[\emph{(}e\emph{)}] For any wall $\ell_i$ there is a map $f\colon T' \to T$ in $\cD(\cT_C)$ such that
\begin{itemize}
\item for any $(b,w) \in \ell_i$, the objects $T', T$ lie in the heart $\cA(b)$,
\item $T$ is $\sigma\_{b,w}$-semistable of class $v$ with $\nu\_{b,w}(T')=\nu\_{b,w}(T)=\,\mathrm{slope}\,(\ell_i)$ constant on the wall $\ell_i$, and
\item $f$ is an injection $T' \hookrightarrow T$ in $\cA(b)$ which strictly destabilises $T$ for $(b,w)$ in one of the two chambers adjacent to the wall $\ell_i$.
\hfill$\square$
\end{itemize} 
	\end{itemize*} 
\end{Prop}

Fix a numerical class $v=(r,d,n)\in\mathbb{Z}_{>0}^{3}$ and $b_0<d/r$, and consider objects of class $v$ in the heart $\mathcal{A}(b_0)$. The following descriptions follow from the identification with $\alpha$-stability for sufficiently small $\alpha>0$; see \cite[Proposition~6.6]{FN} and \cite[Remark~1.4(ii)]{RP94}. The wall description extends this identification from $b_0<0$ to every $b_0<d/r$.

\begin{Lem}\label{lem-LVL}
    An object $T \in \cA(b_0)$ of class $v$ is $\sigma_{b_0,w}$-stable for all $w\gg 0$ if and only if it is represented by a coherent system $ T=[\mathcal{O}_C\otimes V\overset{\ \varphi}{\to}E]\in\T$ satisfying the following conditions: 
\begin{enumerate}
    \item $E$ is a slope-semistable vector bundle,
    \item the induced map $V\to H^0(E)$ is injective,
    \item for every nonzero proper subobject $T'\subset T$ in $\T$ of positive rank, 
\[
  \frac{\mathbf{d}(T')}{\mathbf{r}(T')}=\frac{d}{r} \qquad \Rightarrow \qquad \frac{\mathbf{n}(T')}{\mathbf{r}(T')}<\,\frac{n}{r}.
  \]
\end{enumerate}
\end{Lem}
\medskip

We denote by $M_\infty(v)$ the coarse moduli space of $\sigma_{b_0,\,w\gg 0}$-stable objects of class $v$ in $\mathcal{A}(b_0)$ where $b_0<d/r$.

\newpage

\section{Butler's conjecture in the high-slope range}\label{section.3}

In this section, we prove Butler's conjecture for classes with
slope greater than $2g-2$. 

\begin{Thm}\label{thm-any rank-high slope}
Let $C$ be a smooth complex projective curve of genus $g\ge8$ such that 
\begin{equation}\label{cond on C}
    \Phi_C(b)<1+\frac b4\qquad(0<b\le g-3).
\end{equation}
Pick triple $(r, d, n) \in (\Z_{>0})^{\oplus 3}$ such that 
\begin{equation}\label{assumption}
    n>r>0 \qquad \text{and} \qquad \frac{d}{r}>2g-2.
\end{equation}
If $M_{\infty}(r,d,n)$ is non-empty, then for a general coherent system $T= [\cO \otimes V \xrightarrow{\varphi} E] \in M_\infty(r,d,n)$, the map $\varphi$ is surjective and its kernel
\[
M_{T} :=\ker(\varphi).
\]
is a stable vector bundle. 
\end{Thm}

Curves on K3 surfaces provide examples of curves satisfying \eqref{cond on C}.

\begin{Lem}
  Let $(X,H)$ be a smooth polarised K3 surface such that $H^2$ divides $H.D$ for every curve class $D$ on $X$, and let $C\in|H|$ be a smooth curve of genus $g \geq 8$. Then $C$ satisfies \eqref{cond on C}.
\end{Lem}
    \begin{proof}
If $0<b<1$, every semistable bundle $F$ of slope $b$ satisfies $H^0(F(-p))=0$ for any $p\in C$. Hence evaluation at $p$ is injective on $H^0(F)$, and $h^0(F)\le\rk(F)$. Passing to the defining upper limit gives
\[
  \Phi_C(b)\le 1<1+\frac{b}{4}.
\]
For $1\le b\le g-3$, the quadratic section bound of \cite[Theorem~1.1 and \S4.2]{feyz-li:clifford-indices}, which holds under the stated divisibility hypothesis, gives
\[
  \Phi_C(b)\le 1+\frac1g+\frac{gb^2}{4(g-1)^2}.
\]
Here the inequality is non-strict because we have passed to the upper limit.
Since $(b-1)(b-g+3)\le 0$, we have
\[
  b^2\le(g-2)b-(g-3).
\]
Consequently,
\[
\begin{aligned}
  \Phi_C(b)-\Bigl(1+\frac{b}{4}\Bigr)
  &\le\frac1g+\frac{gb^2}{4(g-1)^2}-\frac{b}{4}\\
  &\le\frac1g-\frac{g(g-3)+b}{4(g-1)^2}\\
  &<\frac1g-\frac{g-3}{4g}<0,
\end{aligned}
\]
where the last inequality uses $g\ge 8$. 
\end{proof}

\medskip
\subsection{Reduction to classes close to the Riemann--Roch line} The following lemma shows that, to prove Theorem~\ref{thm-any rank-high slope}, it suffices to find a single coherent system with the desired properties. The subsequent proposition then reduces the problem to classes close to the Riemann--Roch line.

\begin{Lem}\label{lem-signle}
  Let $C$ be a smooth complex projective curve of genus $g\ge 2$, and suppose that $n>r>0$ and $d/r>2g-2$. Then $M_\infty(v)$ for $v= (r, d, n)$ is nonempty if and only if $n\le d-r(g-1)$. Whenever it is nonempty, it is smooth and irreducible of dimension
  \[
    1- \chi(v, v) = r^2(g-1)+1+n\bigl(d-r(g-1)-n\bigr).
  \]
  Moreover, if $M_\infty(v)$ contains one coherent system with surjective evaluation map and slope-stable kernel, then a general system in $M_\infty(v)$ has both properties.
\end{Lem}

\begin{proof}
By Lemma~\ref{lem-LVL}, every system in $M_\infty(r,d,n)$ has slope-semistable underlying bundle $E$ and determines an $n$-dimensional subspace of $H^0(C,E)$. Since $\mu(E)>2g-2$, Serre duality and Riemann--Roch give
\[
  H^1(C,E)=0,\qquad h^0(C,E)=d-r(g-1).
\]
Thus $n\le d-r(g-1)$ is necessary. Conversely, when this inequality holds, choose a slope-stable bundle $E$ of rank $r$ and degree $d$, then $H^1(C,E)=0$. Every $n$-dimensional subspace of $H^0(C,E)$ then defines a large volume stable system by Lemma~\ref{lem-LVL}.

By cohomology and base change, these spaces of sections form a vector bundle over the irreducible stack of slope-semistable bundles of rank $r$ and degree $d$ \cite[Corollary~A.5 and the following paragraph]{Hoffmann2010}. Its relative Grassmannian of $n$-planes is therefore irreducible. Large volume stable systems form a nonempty open substack of this Grassmannian, so their moduli space $M_\infty(r,d,n)$ is irreducible. It is smooth of the stated dimension by the Petri criterion \cite[Proposition~3.10]{BradlowGaraPradaMunozNewstead2003}, since the source of the Petri map,
\[
  V\otimes H^0(C,E^\vee\otimes\omega_C),
\]
vanishes. Finally, surjectivity of the evaluation map is an open condition in families. On this open locus, the kernels form a family of vector bundles, so slope stability of the kernel is also open. A single generated system with slope-stable kernel therefore gives a nonempty open subset of $M_\infty(r,d,n)$, which is dense by irreducibility.
\end{proof}

\bigskip

The first step to prove Theorem~\ref{thm-any rank-high slope} is reducing to classes in a narrow strip below the
Riemann--Roch line. 

\begin{Prop}[Increasing the degree by the kernel rank]\label{prop:degree-shift-gen}
  Let $C$ be a smooth complex projective curve of genus $g\ge 2$. Suppose there is an object
  \[
    T=\bigl[\mathcal{O}_C\otimes V\xrightarrow{\ \varphi\ }E\bigr] \in \cT_C
  \]
  of class $(r,d,n)$ with $n>r>0$ such that $E$ is a vector bundle, $\varphi$ is surjective, $H^1(C,E)=0$, and $\ker\varphi$ is slope-stable.
  Then there exists a coherent system
  \[
    T'=\bigl[\mathcal{O}_C\otimes V\xrightarrow{\ \varphi'\ }E'\bigr]
  \]
  of class $(r,\,d+n-r,\,n)$ such that $\varphi'$ is surjective, both $E'$ and $\ker\varphi'$ are slope-stable, and $H^1(C,E')=0$. 
\end{Prop}

\begin{proof}
Put $M=\ker\varphi$ and $m=n-r$. Fix a point $p\in C$ and consider
\[
  Q=(\mathcal{O}_C\otimes V)/M(-p).
\]
Then it fits in the following exact sequence
\[
  0\longrightarrow\mathcal{O}_p^{\oplus m}\longrightarrow Q\longrightarrow E\longrightarrow 0.
\]
Thus $Q$ has rank $r$ and degree $d+m$, and $H^1(C,Q)=0$. Its evaluation kernel $M(-p)$ is stable. 

Let $\mathcal{U}$ be the locus in
\[
  \operatorname{Quot}(\mathcal{O}_C\otimes V;\,r,\,d+m)
\]
consisting of quotients with vanishing $H^1$ and slope-stable kernel. This locus is open by semicontinuity \cite[III, Theorem~12.8]{Hartshorne} and the openness of stability \cite[Proposition~2.3.1]{HL}. The quotient $Q$ constructed above belongs to $\mathcal{U}$, so $\mathcal{U}\ne\varnothing$.

Consider the forgetful morphism
\[
  \mathcal{U}\longrightarrow\operatorname{Coh}_{r,d+m}^{H^1=0},
\]
where the target is the stack of coherent sheaves of rank $r$ and degree $d+m$ with vanishing $H^1$. We claim this morphism is smooth. Indeed, by cohomology and base change \cite[III, Theorem~12.11]{Hartshorne}, the stack of morphisms $\mathcal{O}_C\otimes V\to F$ is a vector bundle over $\operatorname{Coh}_{r,d+m}^{H^1=0}$, with fibre $\operatorname{Hom}\bigl(V,H^0(C,F)\bigr)$. Surjections form an open substack of this bundle, and $\mathcal{U}$ is open therein. Thus the forgetful morphism is the restriction of the vector bundle projection to an open substack, and hence is smooth.

Consequently, the image of $\mathcal{U}$ is a nonempty open substack of $\operatorname{Coh}_{r,d+m}$. By~\cite[Corollary~A.5 and the following paragraph]{Hoffmann2010}, $\operatorname{Coh}_{r,d+m}$ is irreducible, and its locus of stable vector bundles is nonempty and open for $g\ge 2$, hence it intersects the image of $\mathcal{U}$. Choose a quotient $\varphi'\colon\mathcal{O}_C\otimes V\twoheadrightarrow E'$ in this intersection. Then $E'$ and $\ker\varphi'$ are slope-stable and $H^1(C,E')=0$. 
\end{proof}

\medskip

\subsection{Controlling the obstruction}
Before starting the wall-crossing analysis for the proof of Theorem~\ref{thm-any rank-high slope}, we state two auxiliary lemmas that will help us control the obstruction space. We first use Serre duality to control the obstruction by relating it to the number of global sections of a suitable vector bundle on $C$.

\begin{Lem}\label{lem: bound on ext2}
Let $F$ be a vector bundle on a smooth projective curve $C$,
and let $W\subset H^0(F)$ be an $n$-dimensional subspace.
Consider the coherent system
\[
T=[\cO\otimes W\xrightarrow{\varphi = \mathrm{ev}|_W}F],
\qquad M=\ker\varphi.
\]
Put
\[
S=\operatorname{im}\!\big(
H^0(F^\vee \otimes \omega_C)\otimes\cO\longrightarrow F^\vee \otimes \omega_C
\big).
\]
Then
\[
\begin{aligned}
\ext^2(T,T)=h^0(M\otimes F^\vee \otimes \omega_C) \le\rk(M)\bigl(h^0(S)-\rk(S)\bigr).
\end{aligned}
\]
\end{Lem}

\begin{proof}
Since $W\subset H^0(F)$, we have $H^0(M)=0$.
Moreover, $\operatorname{im}\varphi$ is torsion-free on a smooth
curve, hence locally free. Thus $M\subset\cO\otimes W$ is a
subbundle, and dualizing its defining sequence gives a surjection
\[
W^\vee\otimes\cO\twoheadrightarrow M^\vee.
\]
In particular, $M^\vee$ is globally generated. By \cite[Proposition~2.2]{FN} and Serre duality,
\[
\begin{aligned}
\Ext^2(T,T)^\vee
&\cong
\ker\!\big(
W\otimes H^0(F^\vee\otimes\omega_C)
\longrightarrow
H^0(F\otimes F^\vee\otimes\omega_C)
\big)\\
&=H^0(M\otimes F^\vee\otimes\omega_C).
\end{aligned}
\]
The last equality follows by tensoring the left exact sequence
\[
0\longrightarrow M\longrightarrow \cO\otimes W
\xrightarrow{\varphi}F
\]
with $F^\vee\otimes\omega_C$ and taking global sections.

We now prove $H^0(M\otimes F^\vee\otimes\omega_C) = H^0(M \otimes S)$. A section of $M\otimes F^\vee\otimes\omega_C$ is a morphism
\[
f:M^\vee\longrightarrow F^\vee\otimes\omega_C.
\]
Precomposing with $W^\vee\otimes\cO\twoheadrightarrow M^\vee$
gives a collection of global sections of
$F^\vee\otimes\omega_C$. Their images lie in $S$ by definition.
Since this precomposition is with a surjection, the image of $f$
itself lies in $S$. Consequently, $f$ factors uniquely through
$S\hookrightarrow F^\vee\otimes\omega_C$.

Conversely, every morphism $M^\vee\to S$ gives such a section.
Hence the natural inclusion induces an isomorphism
\[
H^0(M\otimes S)
=\Hom(M^\vee,S) \xhookrightarrow{\ \simeq\ }
\Hom(M^\vee,F^\vee\otimes\omega_C)=H^0(M\otimes F^\vee\otimes\omega_C).
\]

It remains to prove the required bound. Write
$m=\rk(M)$ and $s=\rk(S)$; the cases $m=0$ or $s=0$ are immediate.
The sheaf $S$ is locally free and globally generated.

Choose $m$ general global sections of $M^\vee$.
Dualizing their evaluation gives an injection
\[
i\colon M\longrightarrow\cO^{\oplus m}
\]
which is an isomorphism away from a finite set $D$.
Since $S$ is globally generated, we can choose $s$ generically
independent sections of $S$ which form a basis at every point
of $D$. Indeed, each basis condition
is a nonempty open subset of $H^0(S)^{\oplus s}$, thus finitely many of those have a nonempty intersection.
We obtain
\[
0\longrightarrow\cO^{\oplus s}
\xrightarrow{j}S\longrightarrow Q\longrightarrow0,
\qquad
\operatorname{Supp}Q\cap D=\varnothing.
\]
In particular, $i$ is an isomorphism near $\operatorname{Supp}Q$, thus it induces $M\otimes Q\simeq Q^{\oplus m}$.

Tensoring this sequence with $M$ and with $\cO^{\oplus m}$ and taking global sections, gives the commutative diagram with exact rows
\[
\begin{tikzcd}[row sep=small]
0 \arrow[r]
& H^0(M^{\oplus s})=0 \arrow[r] \arrow[d]
& H^0(M\otimes S)
    \arrow[r, hook]
    \arrow[d, hook]
& H^0(M\otimes Q)
    \arrow[d, "\simeq"] \\
0 \arrow[r]
& H^0(\mathcal O_C^{\oplus s})^{\oplus m} \arrow[r]
& H^0(S)^{\oplus m} \arrow[r]
& H^0(Q)^{\oplus m}.
\end{tikzcd}
\]
and the injection
\[
H^0(M\otimes S)\hookrightarrow
\left(
H^0(S)/j_*H^0(\cO^{\oplus s})
\right)^{\oplus m}.
\]
Finally, $j_*$ is injective and
$h^0(\cO^{\oplus s})=s$. Taking dimensions gives
\[
h^0(M\otimes S)
\le m\bigl(h^0(S)-s\bigr),
\]
as required.
\end{proof}

The next lemma allows us to use the function $\Phi_C$ to bound the number of global sections of an arbitrary vector bundle on $C$.

\bigskip
\begin{Lem}\label{lem:bound on sections}
Let $F$ be a nonzero vector bundle on a smooth projective
curve $C$, with $\mu^{-}(F)\ge0$.
If $\deg F>0$, then
\[
h^0(F)-\rk(F)
\le
\left(
\sup_{\substack{\mu^{-}(F)\le u\le\mu^{+}(F)\\ u>0}}
\frac{\Phi_C(u)-1}{u}
\right)\deg F.
\]
If $\deg F=0$, then $h^0(F)\le\rk(F)$.
\end{Lem}

\begin{proof}
Let $F_i$ be the semistable factors of the HN filtration
of $F$. Subadditivity of $h^0$ and additivity of rank give
\[
h^0(F)-\rk(F)
\le
\sum_i\bigl(h^0(F_i)-\rk(F_i)\bigr).
\]
For each positive-slope factor,
\[
\begin{aligned}
h^0(F_i)-\rk(F_i)
&\le \rk(F_i)\bigl(\Phi_C(\mu(F_i))-1\bigr)\\
&=\deg(F_i)\,
\frac{\Phi_C(\mu(F_i))-1}{\mu(F_i)}.
\end{aligned}
\]
A slope-zero factor satisfies $h^0(F_i)\le\rk(F_i)$,
so its contribution is nonpositive.
Since all degrees $\deg(F_i)$ are nonnegative, replacing
each quotient by the displayed supremum and summing
proves the inequality.
\end{proof}

The following lemma gives a piecewise linear upper bound for \(\Phi_C\). Its Clifford-index hypothesis follows from \eqref{cond on C}: in genus \(g\ge8\), a curve of Clifford index at most one has a pencil of degree \(b=2\) or \(3\), which gives \(\Phi_C(b)\ge2>1+b/4\).

\begin{Lem}\label{lem:gamma}
Let \(C\) be a smooth complex projective curve of genus \(g\ge8\) with \(\operatorname{Cliff}(C)\ge2\). 
Define, for \(b\ge0\),
\[
\gamma(b) :=
\max\left\{
1,\;
1+\frac{b-1}{g},\;
\frac{(g-3)b+4}{2(g-2)},\;
\frac b2,\;
\frac{3b}{4}+1-\frac{g-1}{2},\;
b+1-g
\right\}.
\]
Then \(\gamma\) is continuous, convex and piecewise linear, and
\[
\Phi_C(b)\le\gamma(b)\qquad(b\ge0).
\]
\end{Lem}
\begin{proof}
    Being the maximum of finitely many affine functions with increasing slopes, $\gamma$ is continuous, convex and piecewise linear. We prove that it bounds $\Phi_C$. The standard small-slope bound \cite[Lemma~3.1]{bradlow-et:moduli-of-coheent-system-of-small-slope} and Mercat's theorem (see \cite[Corollary~3.4 and its proof]{LangeNewstead}) give
\[
  \Phi_C(b)\le
  \begin{cases}
    1+\dfrac{b-1}{g}, & 0<b<2,\\[6pt]
    \max\Bigl\{1+\dfrac{b-1}{g-2},\ \dfrac{b}{2}\Bigr\}, & 2\le b\le 4,\\[6pt]
    \dfrac{b}{2}, & 4\le b\le 2g-6.
  \end{cases}
\]
These bounds pass to the upper-semicontinuous envelope defining $\Phi_C$.

For $2\le b\le 4$, we have
\[
\begin{aligned}
  \frac{(g-3)b+4}{2(g-2)}-\Bigl(1+\frac{b-1}{g-2}\Bigr)&=\frac{(g-5)(b-2)}{2(g-2)}\ge 0,\\
  \frac{(g-3)b+4}{2(g-2)}-\frac{b}{2}&=\frac{4-b}{2(g-2)}\ge 0.
\end{aligned}
\]
Thus the displayed bounds imply
\[
  \Phi_C(b)\le\gamma(b)\qquad(0<b\le 2g-6).
\]
At $b=0$, we have $\Phi_C(0)=1=\gamma(0)$.

The same estimates also give
\[
  \Phi_C(b)\le 1+\frac{b}{4}\qquad(0\le b\le 4).
\]
Indeed, this is immediate for $0<b<2$, while for $2\le b\le 4$,
\[
  \frac{(g-3)b+4}{2(g-2)}-\Bigl(1+\frac{b}{4}\Bigr)=\frac{(g-4)(b-4)}{4(g-2)}\le 0.
\]
Consequently, Serre duality yields, for $2g-6\le b\le 2g-2$,
\[
\begin{aligned}
  \Phi_C(b)&=b+1-g+\Phi_C(2g-2-b)\\
  &\le b+1-g+1+\frac{2g-2-b}{4}\\
  &=\frac{3b}{4}+1-\frac{g-1}{2}\le\gamma(b).
\end{aligned}
\]
Finally, for $b>2g-2$, Riemann--Roch gives
\[
  \Phi_C(b)=b+1-g\le\gamma(b).
\]
This proves the assertion for every $b\ge 0$.
\end{proof}

\subsection{Wall-crossing analysis}\label{sec:wall-crossing} From now on, we fix the class $v= (r, d, n)$ as in Theorem \ref{thm-any rank-high slope}. Suppose that $M_{\infty}(v) \neq \emptyset$, i.e. 
$$
n \leq d+r(1-g). 
$$
Fix
\[
-\frac{1}{n-r}<b_0<0,
\]
We know at the large volume limit $M_{b_0, w \gg 0}(v) = M_{\infty}(v)$. We move down this vertical line from the large volume limit to the Brill--Noether wall $\ell_f$, which is induced by $[\cO\xrightarrow{\ \mathrm{id}\ }\cO]$ and is given by
\begin{equation}\label{final wall}
  \ell_f\colon\qquad w=sb+1,\qquad\text{where }\  s:=\frac{n-r}{d}.
\end{equation}
Along this path, objects in $M_\infty(v)$ meet only finitely many walls. The key point is that if a coherent system remains semistable up to the Brill--Noether wall, then it is generated with semistable kernel.

\begin{Lem}\label{lem-bn}
  Let $T=[\cO\otimes V\xrightarrow{\ \varphi\ }E]\in M_\infty(r,d,n)$. If $T$ is $\sigma_{b_0,w}$-semistable for all $w>sb_0+1$, then $\varphi$ is surjective and $M_T=\ker\varphi$ is slope-semistable.
\end{Lem}

\begin{proof}
Put \(w_0=1+sb_0>0\). By continuity, \(T\) is semistable at \((b_0,w_0)\). Since \(b_0<0\) and the BN wall has positive slope \(s\), the evaluation map \(\varphi\) is surjective by \cite[Lemma~6.5]{FN}.

At this point of the wall, the stable object \([\mathcal O_C\xrightarrow{\mathrm{id}}\mathcal O_C]\) has the same phase as \(T\). In the abelian category of that phase, consider the natural map
\[
[\mathcal O_C\xrightarrow{\mathrm{id}}\mathcal O_C]\otimes V\longrightarrow T.
\]
Its source is a direct sum of copies of one simple object, so any nonzero kernel would contain a copy of that object. But applying \(\operatorname{Hom}([\mathcal O_C\xrightarrow{\mathrm{id}}\mathcal O_C],-)\) to the map gives the identity \(V\to V\). The map is therefore injective. Since \(\varphi\) is surjective, its quotient is \([0\to M_T][1]\), which is consequently semistable at \((b_0,w_0)\).

Suppose \(M_T\) were slope-unstable, and let \(Q\) be its last HN quotient. Then
\[
\mu(Q)<\mu(M_T)=-\frac{d}{n-r}<b_0.
\]
Since \([0\to M_T][1]\in\mathcal A(b_0)\), we have that every HN slope of \(M_T\) is at most \(b_0\). The kernel of \(M_T\twoheadrightarrow Q\) consists of the preceding HN factors, so shifting this exact sequence gives a quotient \([0\to M_T][1]\twoheadrightarrow[0\to Q][1]\) in \(\mathcal A(b_0)\). Its tilt slope satisfies
\[
\nu_{b_0,w_0}\bigl([0\to Q][1]\bigr)
=\frac{w_0}{b_0-\mu(Q)}
<\frac{w_0}{b_0-\mu(M_T)}
=\nu_{b_0,w_0}\bigl([0\to M_T][1]\bigr).
\]
This contradicts the $\sigma_{b_0, w_0}$-semistability of $[0\to M_T][1]$, since every nonzero quotient of a semistable object has tilt slope at least that of the object. Hence \(M_T\) is slope-semistable.
\end{proof}

Fix a wall $\mathcal W$ for the class $v$ encountered as $w$ decreases along the vertical line $b=b_0$, up to and including the Brill--Noether wall $\ell_f$. Write $(b_0,w_0)$ for the point where this line meets $\mathcal W$. For $T\in M_\infty(v)$ that remains semistable down to $\mathcal W$, a \emph{destabilizing sequence} for $T$ on $\mathcal W$ is a short exact sequence
\begin{equation}\label{eq:wall-extension}
  0\longrightarrow T_1\xrightarrow{\ i\ }T\xrightarrow{\ p\ }T_2\longrightarrow 0
\end{equation}
in $\mathcal A(b_0)$ such that $T_1$ and $T_2$ are nonzero $\sigma_{b_0,w_0}$-semistable objects of the same phase as $T$.


By taking cohomology of ~\eqref{eq:wall-extension} with respect to $\T$ we get that $T_1\in \T$. Moreover, $b_0<0$ implies $i^*\mathcal{H}^{-1}(T_2)=0$. Then we have the long exact sequence in $\cT_C$:
\begin{equation}\label{long}
0\rightarrow
\left[\begin{array}{c}0\\ \big\downarrow\\ j^{!}\cH^{-1}(T_2)\end{array}\right]
\rightarrow
T_1=\left[\begin{array}{c}\cO\otimes W\\ \big\downarrow\rlap{$\scriptstyle\psi$}\\ F\end{array}\right]
\rightarrow
T=\left[\begin{array}{c}\cO\otimes V\\ \big\downarrow\rlap{$\scriptstyle\varphi$}\\ E\end{array}\right]
\rightarrow\cH^0(T_2)\rightarrow0 .
\end{equation}
This implies $W\subset V$ and $F$ is a vector bundle.

\begin{Lem}\label{lem:wall-dimension}
  If every destabilizing sequence \eqref{eq:wall-extension} as above satisfies
  \begin{equation}\label{eq:assum}
    \operatorname{Ext}^2(T_2,T)=0
    \qquad\text{and}\qquad
    -\chi(T_1,T_2)>\operatorname{ext}^2(T_1,T_1)+\operatorname{ext}^2(T_2,T_2),
  \end{equation}
  then the locus of systems in $M_\infty(v)$ whose first destabilizing wall is $\mathcal{W}$ has dimension strictly smaller than $\dim M_\infty(v)$.
\end{Lem}

\begin{proof}
Since $b_0<0$, every object of $\mathbb F^{b_0}$ has the form $[0\to A]$ (a triple with $V\neq0$ has a subobject of slope $\ge0>b_0$). By \cite[Proposition~2.2]{FN}, $\Ext^k$ vanishes for $k>2$ on $\cT_C$, and the exact sequence there with $V_1=0$ gives $\Ext^2([0\to A],T')=0$ for all $T'\in\cT_C$. Hence 
\begin{equation}\label{hd}
    \Ext^k(X,Y)=0 \qquad \text{for all} \ X,Y\in \cA(b_0) \ \text{and} \ k \geq 3.
\end{equation}
Since $T$ is large-volume stable, $\Hom(T,T)=\mathbb C$, so every endomorphism of $T$ preserves $T_1$. Moreover, by Lemma \ref{lem-signle}, we know $\Ext^2(T, T) = 0$.

We have $T_1\in\cT_C$, and \eqref{eq:wall-extension} is equivalent to the morphism $i\colon T_1\to T$ of $\cT_C$, whose kernel $\cH^{-1}(T_2)\in\mathbb F^{b_0}$ may be nonzero. Let $\operatorname{Def}(T_1\subset T)$ denote the space of
first-order deformations of this morphism, equivalently of the exact sequence \eqref{eq:wall-extension}.
Applying the filtered derived-Hom construction of
\cite[\S 1.5, Proposition~1.3 and the subsequent remark]{DL}
to the two-step filtration $0\subset T_1\subset T$, which uses only the triangle \eqref{eq:wall-extension} and not that the filtration is by subsheaves, gives
the following long exact sequence; see also
\cite[\S 2, equation~(4)]{Walter}:
\begin{equation}\label{eq.long}
\begin{aligned}
0\longrightarrow\Hom(T_2,T_1)
&\longrightarrow\Hom(T,T)=\mathbb C\\
&\longrightarrow
\Hom(T_1,T_1)\oplus\Hom(T_2,T_2)\\
&\longrightarrow\Ext^1(T_2,T_1)
\longrightarrow\operatorname{Def}(T_1\subset T)\\
&\longrightarrow
\Ext^1(T_1,T_1)\oplus\Ext^1(T_2,T_2)
\xrightarrow{\delta}\Ext^2(T_2,T_1)
\longrightarrow\cdots.
\end{aligned}
\end{equation}
Here $\delta$ is the canonical connecting homomorphism
associated with \eqref{eq:wall-extension}. Its second
component is, up to sign, the boundary map in
\[
\Ext^1(T_2,T_2)\longrightarrow
\Ext^2(T_2,T_1)\longrightarrow
\Ext^2(T_2,T) \overset{\eqref{eq:assum}}{=}0.
\]
Thus $\delta$ is surjective.

Furthermore, $\Hom(T_2,T_1)=0$. Indeed, any morphism
$f\colon T_2\to T_1$ induces the square-zero endomorphism $ifp$
of $T$, where $i\colon T_1\to T$ and $p\colon T\to T_2$ are the maps in \eqref{eq:wall-extension}. Since $T$ is simple, this endomorphism vanishes,
and hence $f=0$, as $i$ is a monomorphism and $p$ an epimorphism in $\mathcal A(b_0)$.

Taking dimensions in \eqref{eq.long} therefore gives
\begin{equation*}
\begin{aligned}
\dim\operatorname{Def}(T_1\subset T)
={}&\ext^1(T_1,T_1)+\ext^1(T_2,T_2)
     +\ext^1(T_2,T_1)\\
&-\hom(T_1,T_1)-\hom(T_2,T_2)+1
     -\ext^2(T_2,T_1)\\
={}&1-\chi(T_1,T_1)-\chi(T_2,T_2)-\chi(T_2,T_1)\\
&+\ext^2(T_1,T_1)+\ext^2(T_2,T_2)\\
={}&1-\chi(T,T)+\chi(T_1,T_2)\\
&+\ext^2(T_1,T_1)+\ext^2(T_2,T_2)\\
 \overset{\eqref{eq:assum}}{<}{}&1-\chi(v,v).
\end{aligned}
\end{equation*}

At a fixed wall, the support property \cite[Lemma~3.7]{FN} leaves only finitely many numerical classes for the wall factors. For each class, the pairs $T_1\subset T$ form a bounded family: $T$ ranges over the finite-type space $M_\infty(v)$, and the sheaf of $T_1$ has fixed rank and degree and $\mu_{\max}$ bounded by that of the sheaf of $T$. The wall locus is therefore covered by the images of finitely many incidence strata parametrising these filtered objects. Their local dimensions are at most $\dim\operatorname{Def}(T_1\subset T)$, and forgetting the filtration cannot increase dimension. On the other hand, since $\Hom(T,T)=\mathbb C$ and $\Ext^2(T,T)=0$, the space $M_\infty(v)$ is smooth of dimension $\ext^1(T,T)=1-\chi(v,v)$ at every point. Hence the wall locus has dimension strictly smaller than $\dim M_\infty(v)$. 
\end{proof}


By analyzing the destabilizing factors in our setting, we in fact obtain the following vanishing of obstructions.

\begin{Lem}\label{lem:ext2}
    We have $\ext^2(T_2, T)=0$ and $\ext^2(T_2, T_2)=0$.
\end{Lem}
\begin{proof}
    Since $E$ is semistable with $\mu(E)>2g-2$, we get $H^1(E)=0$, thus 
    \begin{equation}\label{f1}
        \ext^2(\mathcal{H}^0(T_2), T)=0
    \end{equation}
    follows from~\cite[Proposition 2.2]{FN}. Since
    $i^*\mathcal{H}^{-1}(T_2)=0$ and $\mu^+(j^!\mathcal{H}^{-1}(T_2))\leq b_0<0$, we get that 
    \begin{equation}\label{f2}
        \Hom^2(\mathcal{H}^{-1}(T_2)[1], T)=\Hom^1(j^!\mathcal{H}^{-1}(T_2), E)=\Hom(E\otimes \omega_C^{-1}, j^!\mathcal{H}^{-1}(T_2))=0,
    \end{equation}
    Thus, combining \eqref{f1} and \eqref{f2} with the exact triangle $
\cH^{-1}(T_2)[1]\rightarrow T_2\rightarrow\cH^0(T_2)
$
gives $\Ext^2(T_2,T)=0$.
    
  It remains to prove the second vanishing. Applying $\Hom(T_2,-)$ to~\eqref{eq:wall-extension} gives the exact sequence
\[
0=\Hom^2(T_2,T)\longrightarrow\Hom^2(T_2,T_2)\longrightarrow\Hom^3(T_2,T_1)=0,
\]
where the first vanishing is the first part of the claim and the last follows from~\eqref{hd}. Hence $\Hom^2(T_2,T_2)=0$, which proves the second claim.
    
\end{proof}


\begin{proof} [Proof of Theorem \ref{thm-any rank-high slope}] 
Write $\cl(T_1)=(a,e,k)$ and equation of the wall $\mathcal{W}$ as $ w=tb+c$, then 
\begin{equation*}
    n=cr+td \qquad \text{and} \qquad  k=ca+te. 
\end{equation*}
Since $\mathcal{W}$ lies above $\ell_f$ \eqref{final wall}, we have  
\begin{equation*}
    c>1,\qquad t<\frac{n-r}{d}.
\end{equation*}
Taking $i^{*}$ from \eqref{long} gives 
\begin{equation*}
    k \leq n.
\end{equation*}
Taking $j^{!}$ from \eqref{long} gives the long exact seqeunce of sheaves on $C$:
\begin{equation*}
    0 \to j^{!}\cH^{-1}(T_2) \to F \to E \to j^{!} \cH^0(T_2) \to 0.  
\end{equation*}
Let $b_{\max} \geq 0$ be the largest $b$-coordinate of an intersection point of the final BN wall $\ell_f$ in~\eqref{final wall} with the graph of the convex function $\gamma$. Moving along the wall $\mathcal{W}$, which lies above $\ell_f$, shows that
\begin{equation}\label{min-max}
\mu^+\bigl(j^{!}\cH^{-1}(T_2)\bigr)\le-\frac{d}{n-r}<b_0<0
\qquad\text{and}\qquad
0<b_{\max}<\mu^-\bigl(j^{!}\cH^0(T_2)\bigr).
\end{equation}
We know $\Pi(T_1)$ lies above $\ell_f$, hence
\begin{equation}\label{eq: thm 3.1 - 1}
\frac{k}{a}>s\frac{e}{a}+1>1 \qquad \text{and} \qquad
    0<  b_{\max} < \mu^-(F) \leq \mu(F)= \frac{e}{a} \leq  \mu^+(F) \leq \frac{d}{r}.
\end{equation}
Moreover, we have 
\begin{equation}\label{eq:ad-re}
    ad-re>0.
\end{equation}
Otherwise for a point $(w_0, b_0)$ on the wall $\mathcal{W}$ the equality $\nu_{b_0, w_0}(T_1)=\nu_{b_0, w_0}(T)$ implies that $\Pi(T_1)=\Pi(T)$. Thus $T_1$ destabilizes $T$ for all stability conditions on the vertical line $b=b_0$, which contradicts with $T$ being stable in the large volume limit. 



The first step of the proof is to show that a general object of $M_\infty(r,d,n)$ remains stable after crossing every wall above $\ell_f$ along the vertical line $b=b_0$; Lemma~\ref{lem-bn} then gives the semistability of the kernel. By Lemmas~\ref{lem:wall-dimension} and~\ref{lem:ext2}, to show that a general object of $M_\infty(r, d, n)$ is not destabilized along the given wall $\mathcal{W}$, it suffices to prove that
\begin{equation}\label{goal}
  -\chi(T_1,T_2)>\operatorname{ext}^2(T_1,T_1).  
\end{equation}
A direct calculation of the Euler form gives
\begin{equation}
\begin{aligned}\label{chi}
-\chi(T_1,T_2)
={}&\left(1-\frac kn\right)
 \bigl[re-an+(k-a)(\chi(E)-n)\bigr]\\
&+\frac{ad-re}{n}
 \left[k(c-1)+t(g-1)(k-a)\right].
\end{aligned}
\end{equation}
To prove \eqref{goal}, we first show that the first bracket in \eqref{chi} is positive, and then use the second summand to control the obstruction.

 \bigskip

\noindent\emph{Step 1. Positivity of the first bracket in \eqref{chi}.} Dividing the first bracket in~\eqref{chi} by $ar>0$, we obtain
\[
\begin{aligned}
\frac{re-an+(k-a)(\chi(E)-n)}{ar}
&=\frac ea-\frac nr
+\left(\frac ka-1\right)\frac{\chi(E)-n}{r}\\
& \overset{\eqref{eq: thm 3.1 - 1}}{\geq} \frac{e}{a} \left( 1+ s\ \frac{\chi(E) -n}{r}\right) - \frac{n}{r}\\
&\overset{\eqref{eq: thm 3.1 - 1}}{>}
b_{\max}\left( 1+ s\ \frac{\chi(E) -n}{r}\right)-\frac nr.
\end{aligned}
\]
Define
\begin{equation*}\label{eq def b_*}
b_*:=\frac{\frac{n}{r}}{1+s\,\frac{\chi(E)-n}{r}} = \frac{n}{(1-s)n+sr(1-g)} >0.
\end{equation*}
Then to prove positivity of the first bracket, it suffices
to show $b_{\max}>b_*$. Since $\gamma$ is convex and $b_{\max}$ is the positive intersection
of $\gamma$ with the BN line $\ell_f$, it is enough to prove
\begin{equation}\label{b-condition}
\gamma(b_*)<1+sb_*.
\end{equation}

Our assumptions imply
\begin{align}\label{c1}
    0 < s= \frac{n-r}{d} \leq \frac{\chi(E)-r}{d} = 1 -g\frac{r}{d} < 1. 
\end{align}

Moreover, 
\begin{equation}\label{c2}
    \frac{n-r}{2} = \frac{sd}{2} > sr(g-1),
\end{equation}
and 
\begin{equation}\label{c3}
    (1-s)n -srg = r(1-s) +s (\chi(E) -n) \geq r(1-s).
\end{equation}
This implies $0<(1-s)n-sr(g-1)<n$ and so $b_*>1$. Then constant branch $1$ of $\gamma$ is dominated by $1+(b_*-1)/g$. It therefore remains to compare the other five affine branches with $1+sb_*$. Three comparisons are immediate:
\begin{equation*}\label{eq:easy}
\begin{aligned}
  \frac{n}{b_*}\Bigl[1+sb_*-\Bigl(1+\frac{b_*-1}{g}\Bigr)\Bigr]
    &=\frac{s(g-1)(n-r)}{g}>0,\\[4pt]
  \frac{n}{b_*}\Bigl[1+sb_*-\frac{b_*}{2}\Bigr]
    &=\frac{n}{2}-sr(g-1) \overset{\eqref{c2}}{>}  \frac{r}{2}>0,\\[4pt]
  \frac{n}{b_*}\Bigl[1+sb_*-(b_*+1-g)\Bigr]
    &=(g-1)\bigl((1-s)n-srg\bigr)\\
    & \overset{\eqref{c3}}{ \geq } r(g-1)(1-s)>0.
\end{aligned}
\end{equation*}

For the branch $\bigl((g-3)b+4\bigr)/\bigl(2(g-2)\bigr)$, we have
\begin{equation*}\label{eq:branch4}
\begin{aligned}
  \frac{2(g-2)n}{b_*}\Bigl[1+sb_*-\frac{(g-3)b_*+4}{2(g-2)}\Bigr]
  & =(g-5+4s)n-2sr(g-1)(g-4)\\
  & \overset{\eqref{c2}}{>} (g-5+4s)(2sr(g-1)+r) -2sr(g-1)(g-4)\\ 
  & =r\Bigl[8(g-1)\Bigl(s-\frac18\Bigr)^2+4s+\frac{7g-39}{8}\Bigr]>0,
\end{aligned}
\end{equation*}
since $g\ge 8$. Finally, for the branch $3b/4+1-(g-1)/2$, we have 
\begin{equation*}
\begin{aligned}
  &\frac{4n}{b_*}\Bigl[1+sb_*-\Bigl(\frac{3b_*}{4}+1-\frac{g-1}{2}\Bigr)\Bigr]\\
  &=\bigl[2g-5-2(g-3)s\bigr]n-2sr(g-1)^2\\
  & \overset{\eqref{c3}}{\geq} \bigl[2g-5-2(g-3)s\bigr] \left( r+ \frac{srg}{1-s}\right) -2sr(g-1)^2\\
  &=\frac{r}{1-s}\Bigl[4(g-1)\Bigl(s-\frac58\Bigr)^2+4s+\frac{7g-55}{16}\Bigr]>0,
\end{aligned}
\end{equation*}
because $2g-5-2(g-3)s>1$ by~\eqref{c1} and $g\ge 8$. This completes the proof of \eqref{b-condition}.

\medskip

\noindent\emph{Step 2. The second summand controls the obstruction.} We now prove
\begin{equation}\label{eq:star}
 \operatorname{ext}^2(T_1,T_1)<\frac{ad-re}{n}\bigl[k(c-1)+t(g-1)(k-a)\bigr].
\end{equation}
Put
\[
  M=\ker(W\otimes\mathcal{O}_C\to F),\qquad
  S=\operatorname{im}\Bigl[H^0(F^\vee\otimes \omega_C)\otimes\mathcal{O}_C\longrightarrow F^\vee\otimes \omega_C\Bigr].
\]
Lemma~\ref{lem: bound on ext2} gives
\begin{equation}\label{eq:obs}
  \operatorname{ext}^2(T_1,T_1)=h^0(M\otimes F^\vee\otimes \omega_C)\le\rk(M)\bigl(h^0(S)-\rk(S)\bigr).
\end{equation}
Since $S$ is globally generated,
\begin{equation}\label{eq:deg}
      h^0(S)-\rk S\le\deg S.
\end{equation}
On the other hand, 
$$
\mu\left(\frac{F^{\vee} \otimes \omega_C}{S}\right) \geq \mu^-(F^{\vee} \otimes \omega_C) = 2g-2 - \mu^+(F) \overset{\eqref{eq: thm 3.1 - 1}}{\geq} 2g-2- \frac{d}{r},
$$
which implies 
\begin{equation}\label{eq:degS}
  \deg S\le\frac{ad-re}{r}-\Bigl(\frac{d}{r}-2g+2\Bigr)\rk S. 
\end{equation}

\medskip

\textbf{Case $t\le 0$.} Combining \eqref{eq:deg} with \eqref{eq:degS} implies that we only need to show 
\begin{equation}\label{final}
    \frac{n}{r} \rk M \leq k(c-1) +t(g-1)(k-a).
\end{equation}
Let $Q=\operatorname{coker}(W\otimes\mathcal{O}_C\to F)$. If $Q \neq 0$, consider the point $(b_0, w_0) \in \mathcal{W}$ along the wall,  we have the surjection $T_1 \twoheadrightarrow  [0 \to Q]$ in $\cA(b_0)$ and so $\sigma_{b_0, w_0}$-semistability of $T_1$ implies that 
$$
\nu_{b_0, w_0}(T_1) = t \leq \nu_{b_0, w_0}([0 \to Q]) = \frac{-w_0\rk Q}{\deg Q -b_0 \rk Q}
$$
Since $w_0=tb_0+c$, this yields $c\rk Q+t\deg Q\le 0$, which also holds trivially when $Q=0$. Using $k=ca+te$, $\rk M=k-a+\rk Q$ and $\deg M=-e+\deg Q$, we obtain
\[
  c\rk M+t\deg M=k(c-1)+c\rk Q+t\deg Q\le k(c-1).
\]
Now $M\subset W\otimes\mathcal{O}_C$ gives $\deg M\le 0$. Hence $t\deg M\ge 0$, and so $c\rk M\le k(c-1)$. Finally, since $d/r>g-1$ and $\rk M\ge k-a$, multiplying by $t\le 0$ gives
\[
  t\,\frac{d}{r}\rk M\le t(g-1)(k-a).
\]
Adding the two bounds gives
\[
  \frac{n}{r}\rk M=\left(c+t\,\frac{d}{r} \right)\rk M\le k(c-1)+t(g-1)(k-a),
\]
and so \eqref{final} follows. If $S=0$, the obstruction vanishes, while \eqref{final} with the positivity from~\eqref{eq:ad-re} makes the right-hand side of \eqref{eq:star} positive.

\medskip
For the remaining cases $t>0$, the evaluation map is surjective by \cite[Lemma~6.5]{FN}. Thus
\begin{equation}\label{eq:Mtgt0}
  \rk M=k-a,\qquad \deg M=-e.
\end{equation}

\medskip
\textbf{Case $t>0$ and $s\ge 1/4$.} By our assumption, 
\[
\Phi_C(b)<1+\frac b4\qquad(0<b\le g-3).
\]
Moreover, the wall lies above this line as $s>t$ and $s \geq \frac{1}{4}$, so
\[
tb+c=1+sb+(s-t)\Bigl(\frac dr-b\Bigr)>1+\frac b4.
\]
This allows the wall to be continued to $b=g-3$. Hence $\mu^-(F)>g-3$ and so
\begin{equation*}\label{hn}
    0\le\mu^-(S)\le\mu^+(S)<g+1.
\end{equation*}
It yields $\deg S\le(g+1)\operatorname{rk}S$, alongside  with \eqref{eq:degS} gives
\begin{equation}\label{deg}
    \deg S\le\frac{(ad-re)(g+1)}{d-r(g-3)} .
\end{equation}
Thus, Lemma~\ref{lem:bound on sections} with convexity of $\gamma$ implies that 
\begin{equation}\label{sum}
    h^0(S)-\rk(S)\leq \frac{\gamma(g+1)-1}{g+1}\deg(S)=\frac{g-1}{2(g+1)}\deg(S)
\end{equation}
because, at the point $b=g+1$, the term $b/2$ is the largest of the affine terms defining $\gamma$.
Therefore, \eqref{sum} together with \eqref{deg} imply 
\begin{align}
 \rk (M) (h^0(S) - \rk S) 
&\le(k-a)(ad-re)\,\frac{g-1}{2(d-r(g-3))}\nonumber \\
&<(k-a)(ad-re)\,\frac{s(g-1)}{n} \label{1}\\
&<\frac{ad-re}{n}\bigl[k(c-1)+t(g-1)(k-a)\bigr] \label{2},
\end{align}
where \eqref{1} follows from
\[
2s(d-r(g-3))-n=s(d-2r(g-1))+r(4s-1)>0,
\]
and \eqref{2} comes from
\[
\begin{aligned}
&k(c-1)+t(g-1)(k-a)-s(g-1)(k-a)\\
&\qquad=(s-t)\left( k\,\frac dr-(g-1)(k-a)\right)>0 .
\end{aligned}
\]
This completes the proof of \eqref{eq:star} in this case as $k>a$ by~\eqref{eq: thm 3.1 - 1} and $ad-re>0$ by~\eqref{eq:ad-re}. If $S=0$, then then obstruction vanishes, while~\eqref{2} makes right-hand side of~\eqref{chi} positive.

\medskip

\textbf{Case $0<s<1/4$ and $t>0$.} We may assume $n/r\le g+1$. Otherwise, choose
\[
d_0\equiv d\pmod{n-r},\qquad n+r(g-1)\le d_0<2n+r(g-2).
\]
Since $n\le\chi(E)=d-r(g-1)$, this choice satisfies $d_0\le d$. Moreover,
\[
\frac{d_0}{r}\ge\frac nr+g-1>2g
\qquad\text{and}\qquad
n\le d_0-r(g-1),
\]
so the reduced class $(r,d_0,n)$ still lies in the high-slope range and satisfies our assumptions. Furthermore,
\[
\frac{n-r}{d_0}>\frac{n/r-1}{2n/r+g-2}>\frac13,
\]
where the last inequality is equivalent to $n/r>g+1$. Thus repeated application of Proposition~\ref{prop:degree-shift-gen} reduces this case to the previously established case.

\medskip

In this case, we first control the HN factors of the bundle $M$.

\textbf{Claim 1.} The first HN subbundle $M_0$ of $M$ has slope $-u_0\le-2$. Assume otherwise, i.e.\ $u_0<2$. Let $W_0\subset W$ be the smallest subspace such that $M_0\subset W_0\otimes\cO$, and let $F_0=\operatorname{im}(W_0\otimes\cO\to F)$; note that $F_0\neq0$, since $W\to H^0(F)$ is injective. Thus we have the injection
\[
T_0:=[W_0\otimes\cO\xrightarrow{\psi_0}F_0]\hookrightarrow T_1
\]
in $\cA(b_0)$ with $M_0\subset\ker\psi_0$. Then $\sigma_{b_0,w_0}$-semistability of $T_1$ gives $\nu_{b_0,w_0}(T_0)\le\nu_{b_0,w_0}(T_1)=t$, which implies
\begin{equation}\label{w0}
  \frac{\dim W_0}{\rk F_0}\le c+t\mu(F_0)\le c+t\,\frac dr=\frac nr\le g+1 .
\end{equation}
The bundle $M_0^\vee$ is semistable of slope $u_0<2$, so by \cite[Lemma~3.1]{bradlow-et:moduli-of-coheent-system-of-small-slope},
\[
\dim W_0\le h^0(M_0^\vee)<\frac{g+1}{g}\rk M_0 .
\]
Note that if $u_0=0$ the inequality still holds as $h^0(M_0^\vee)\leq \rk(M_0)$. Hence $\rk F_0\le\dim W_0-\rk M_0<\dim W_0/(g+1)$, contradicting~\eqref{w0}.

\textbf{Claim 2.} Let $M_\ell$ be the last HN quotient of $M$, of slope $-u_\ell$. Then $tu_\ell\le c$. Indeed, by Claim~1 we have $\mu^+(M)\le-2<b_0$, so $[0\to M]\in\mathbb F^{b_0}$ and the exact triangle
\[
[W\otimes\cO\xrightarrow{\mathrm{id}}W\otimes\cO]\to T_1\to[0\to M][1]
\]
is a short exact sequence in $\cA(b_0)$. Thus $[0\to M_\ell][1]$ is a quotient of $T_1$ in $\cA(b_0)$, and
\[
\nu_{b_0,w_0}([0\to M_\ell][1])=\frac{tb_0+c}{u_\ell+b_0}\ge\nu_{b_0,w_0}(T_1)=t .
\]
Since $u_\ell+b_0>0$, this gives $tu_\ell\le c$.
\medskip

Let $M_i$ and $F_j$ be the HN factors of $M$ and $F^\vee \otimes \omega_C$, with slopes $-u_i$ and $v_j$ respectively. The above two claims plus \eqref{eq: thm 3.1 - 1} give
\begin{equation}\label{4}
u_i\ge2,\qquad tu_i\le c,\qquad 2g-2-\frac dr\le v_j<2g-2 .
\end{equation}
We claim that
\begin{equation}\label{eq:Phi(v-u)}
\frac nr\,\Phi_C(v-u)\le\Bigl(v+\frac dr-2g+2\Bigr)\big(c+t(g-1-u)\big),
\end{equation}
for $u,v$ having restrictions as in~\eqref{4}, and it is strict whenever $v+d/r-2g+2>0$.

\medskip

For $v<u$ this is immediate: the left side in~\eqref{eq:Phi(v-u)} is zero and, since $tu\le c$,
\[
c+t(g-1-u)\ge t(g-1)>0.
\]
Suppose now that $u\le v$. We know
\[
\begin{aligned}
&\frac{v+d/r-2g+2}{c+td/r}-\frac{v}{c+2t(g-1)}\\
&\qquad=\frac{(d/r-2g+2)\bigl(c+t(2g-2-v)\bigr)}{(c+td/r)(c+2t(g-1))}>0 .
\end{aligned}
\]
Using $n/r=c+td/r$ and $\Phi_C\le\gamma$, it therefore suffices to prove
\begin{equation}\label{cl}
\bigl(c+2t(g-1)\bigr)\gamma(v-u)<v\bigl(c+t(g-1-u)\bigr).    
\end{equation}
The right-hand side minus the left-hand side is concave in $v$, so it is enough to check $v=u$ and $v=2g-2$. Throughout, we use
\begin{equation}\label{c4t}
    c>4t,
\end{equation}
which follows from $c>1$ and $0<t<s<1/4$.

At $v=u$, using~\eqref{4} and~\eqref{c4t}, we get that the difference is
\[
(u-2)\bigl(c-tu+t(g-3)\bigr)+c-4t>0.
\]

At $v=2g-2$, we distinguish two short cases.

If $2\le u\le4$, the fifth branch of $\gamma$ gives
\[
\gamma(2g-2-u)=g-\frac{3u}{4}.
\]
Consequently, the endpoint difference is affine in $u$. Its values at $u=2$ and $u=4$ are, respectively,
\[
c\Bigl(g-\frac12\Bigr)-3t(g-1)>0,
\qquad
c(g+1)-4t(g-1)>0,
\]
by $c>4t$ from \eqref{c4t}. Hence it is positive throughout this interval.

If $4\le u\le 2g-2$, then one can easily check from Lemma~\ref{lem:gamma}, that 
\[
\gamma(2g-2-u)\le\max\Bigl\{2,\ g-1-\frac u2\Bigr\}.
\]
Therefore, since $g\geq 8$, we get
\begin{equation*}
    (2g-2+u)\,\gamma(2g-2-u)\le\max\Bigl\{2(2g-2+u),\ 2(g-1)^2-\frac{u^2}{2}\Bigr\}<2(g-1)^2.
\end{equation*}
Thus
\[
\frac{\gamma(2g-2-u)}{2g-2}<\frac{g-1}{2g-2+u}\le\frac{c+t(g-1-u)}{c+2t(g-1)},
\]
where the last inequality follows from $tu\le c$. 
This proves the second endpoint inequality, and hence \eqref{cl}.

\medskip

Finally,~\eqref{eq:obs} together with tensoring and summing the HN filtrations for $M$ and $F^\vee\otimes \omega_C$ gives
\[
\begin{aligned}
\operatorname{ext}^2(T_1,T_1)
&=h^0(M\otimes F^\vee\otimes\omega_C)\\
&\le\sum_{i,j}\operatorname{rk}M_i\,\operatorname{rk}F_j\,\Phi_C(v_j-u_i)\\
&\overset{\eqref{eq:Phi(v-u)}}{<}\frac rn\Bigl[\sum_j\operatorname{rk}F_j\Bigl(v_j+\frac dr-2g+2\Bigr)\Bigr]
\Bigl[\sum_i\operatorname{rk}M_i\bigl(c+t(g-1)-tu_i\bigr)\Bigr]\\
&=\frac{ad-re}{n}\bigl[(c+t(g-1))\operatorname{rk}M+t\deg M\bigr]\\
&=\frac{ad-re}{n}\bigl[k(c-1)+t(g-1)(k-a)\bigr].
\end{aligned}
\]
The strict inequality follows because the weights $v_j+d/r-2g+2$ are nonnegative with positive weighted sum $(ad-re)/r$, see~\eqref{eq:ad-re}, and $\operatorname{rk}M>0$. This proves the required obstruction bound \eqref{eq:star}.

\medskip

\emph{Step 3. Stability of the general kernel.} By Steps~1 and~2, a general system \(T=[\mathcal O_C\otimes V\xrightarrow{\varphi}E]\in M_{\infty}(r,d,n)\) reaches the final wall \(\ell_f\), by Lemma~\ref{lem-bn} the map \(\varphi\) is surjective, and \(M_T=\ker\varphi\) is semistable. Suppose that \(M_T\) admits a nonzero proper saturated subbundle
\[
M\subset M_T,\qquad
\mu(M)=-\frac{d}{n-r}.
\]
Define
\[
W=\operatorname{im}\!\left(V^\vee\longrightarrow H^0(M^\vee)\right)^\vee
\subset V,
\qquad
F=(\cO \otimes W)/M, 
\]
thus $W$ is the smallest subspace with $M\subset  \cO \otimes W$, which implies that $F$ is a generated vector bundle with $H^0(F^\vee)=0$. We obtain the commutative diagram in $\mathcal{D}(C)$:
\[
\begin{tikzcd}
  M \arrow[r, hook] \arrow[d]
    & M_T \arrow[r, two heads] \arrow[d]
    & M_T/M \arrow[d] \\
  \mathcal{O}_C\otimes W \arrow[r, hook] \arrow[d]
    & \mathcal{O}_C\otimes V \arrow[r, two heads] \arrow[d, "\varphi"]
    & \mathcal{O}_C\otimes(V/W) \arrow[d] \\
  F \arrow[r, "f"]
    & E \arrow[r]
    & \ker(f)[1] \oplus \operatorname{coker}(f).
\end{tikzcd}
\]
This defines a subobject \(T_1 := [\cO \otimes W \to F]\hookrightarrow T\) in $\mathcal{A}(b_0)$. Indeed, $\ker(f)$ is a subsheaf of the semistable bundle \(M_T/M\), hence \(\mu^+(\ker(f))\leq -d/(n-r)<b_0\), and $\operatorname{coker}(f)$ is a quotient of slope-semistable bundle $E$, thus $\mu^{-}(\text{coker}(f))\geq d/r>0>b_0$. Write $\cl(T_1)=(a,e,k)$, then from $\mu(M)=\mu(M_T)$ it follows
\[
\frac e{k-a}=\frac{d}{n-r},
\]
hence \(T_1\) has the same phase as \(T\) on the BN wall $\ell_f$ and is therefore semistable there.

We now recover the slope bounds needed in Steps~1 and~2. First, 
\begin{equation*}
    \frac{e}{a}= \mu(F) \leq \mu^+(F) \leq \frac{d}{r}
\end{equation*}
because \(F\) is an extension of a subbundle of $E$ and $M_T/M$. Hence \(ad-re\ge0\). Equality would give
\[
\frac ea=\frac dr,\qquad \frac ka=1+s\,\frac dr=\frac nr.
\]
The proper subobject \(T_1\subset T\) in the fixed heart \(\mathcal A(b_0)\) would then have the same phase as \(T\) for all \(w\gg0\), contradicting large-volume stability. Thus \(ad-re>0\).

\medskip

Next, let \(F'\) be the last HN quotient of \(F\). Its slope is positive, because a generated quotient of slope zero is trivial and \(H^0(F^\vee)=0\). The induced coherent system on \(F'\), whose space of sections is the image of \(W\) in \(H^0(F')\), is a quotient of \(T_1\) in \(\mathcal A(b_0)\). Indeed, the kernel of \(F\to F'\) is the preceding HN truncation, whose slopes are all greater than \(b_0\). Semistability of \(T_1\) on \(w=1+sb\) therefore gives
\begin{equation}\label{eq:ker-stable}
    1+s\mu^-(F)\ \le\ \frac{h^0(F')}{\operatorname{rk}F'}\ \le\ \Phi_C\bigl(\mu^-(F)\bigr).
\end{equation}

By Step~1, \(\gamma(b_*)<1+sb_*\). Since \(\gamma\) is convex and \(\gamma(0)=1\), we have \(\gamma(b)<1+sb\) for every \(0<b\le b_*\). Comparing with~\eqref{eq:ker-stable} gives \(\mu^-(F)>b_*\). Hence the first bracket of the Euler formula in~\eqref{chi} is positive, exactly as in Step~1. In fact, \(k-a=se\) gives
\[
re-an+(k-a)(\chi(E)-n)
=a\bigl[r+s(\chi(E)-n)\bigr]\Bigl(\frac ea-b_*\Bigr)>0.
\]

For the second summand in~\eqref{chi}, we use the two positive-slope cases of Step~2.
\begin{itemize}
\item If \(s\ge1/4\), then~\eqref{eq:ker-stable} and the hypothesis \(\Phi_C(b)<1+b/4\) for \(0<b\le g-3\) imply \(\mu^-(F)>g-3\). The degree and section estimates of Step~2 therefore apply to this \(F\) and give
\[
\operatorname{ext}^2(T_1,T_1)
\le(k-a)(ad-re)\,\frac{g-1}{2(d-r(g-3))}
<\frac{ad-re}{n}\,s(g-1)(k-a).
\]
The strict inequality is precisely inequality~\eqref{1} of Step~2, since \(k-a>0\) and \(ad-re>0\). This also covers \(s=1/4\).

\item If \(s<1/4\), the evaluation kernel of \(T_1\) is \(M\), which is already semistable of slope \(-1/s<-4\). Thus the kernel bounds in~\eqref{4} hold directly. At \(c=1\) and \(t=s\), we have \(c>4t\) and \(t\cdot(1/s)=c\). Also \(0<\mu^-(F)\le\mu^+(F)\le d/r\), as required in~\eqref{4}. The factorwise estimate and HN summation at the end of Step~2 consequently give the same strict bound
\[
\operatorname{ext}^2(T_1,T_1)<\frac{ad-re}{n}\,s(g-1)(k-a).
\]
\end{itemize}
Since \(k\le n\), substituting \(c=1\) and \(t=s\) into the Euler formula now gives
\[
-\chi(T_1,T_2)\ \ge\ \frac{ad-re}{n}\,s(g-1)(k-a)\ >\ \operatorname{ext}^2(T_1,T_1).
\]

It remains to apply the deformation count to these refinements. The same argument as in the proof of Lemma~\ref{lem:ext2} shows that the quotient \(T_2=T/T_1\) in \(\cA(b_0)\) satisfies \(\ext^2(T_2,T)=0\) and \(\ext^2(T_2,T_2)=0\). Hence Lemma~\ref{lem:wall-dimension}, together with the inequality \(-\chi(T_1,T_2)>\ext^2(T_1,T_1)\) proved above, shows that the locus of systems admitting such a refinement has positive codimension in \(M_\infty(r,d,n)\) for each numerical type. There are only finitely many types, since \(\operatorname{rk}M\) determines \(\deg M\) and \(\dim W\le n\). A general \(T\) therefore admits no such \(M\), and \(M_T\) is slope-stable.

\end{proof}

\bigskip

\section{Rank one case}\label{section.4}
In this section, we prove Butler's conjecture in rank one and degree at most \(2g-2\) for a general smooth curve in the primitive linear system \(|H|\) on a K3 surface \(X\) with \(\operatorname{Pic}(X)=\mathbb{Z}H\).

We first recall the classical results on the Brill--Noether theory of line bundles on these curves.

\begin{Lem}\label{lem:K3-BN}
Let \(X\) be a smooth K3 surface with \(\operatorname{Pic}(X)=\mathbb{Z}H\), and let \(C\in|H|\) be a general smooth curve of genus \(g\ge 2\). Fix integers \(d\ge 0\) and \(n\ge 1\), and put
\[
  \rho=\rho(g,n-1,d)=g-n(g-d+n-1).
\]
Then:
\begin{enumerate}
  \item If \(\rho<0\), the moduli space \(M_\infty(1,d,n)\) is empty.
  \item If \(\rho=0\), it is nonempty, finite and reduced.
  \item If \(\rho>0\), it is nonempty, smooth and irreducible of dimension \(\rho\).
\end{enumerate}
\end{Lem}
\begin{proof}
Let
\[
  T=[\mathcal{O}_C\otimes V\xrightarrow{\ \varphi\ }L]\in M_\infty(1,d,n).
\]
By \cite[Proposition~2.2]{FN} and Serre duality 
\[
  \operatorname{Ext}^2(T,T)^\vee\cong
  \ker\!\big(V\otimes H^0(L^\vee\otimes K_C)\xrightarrow{\ \mu_T\ }H^0(K_C)\big),
\]
where $\mu_T(v\otimes s)=H^0(\varphi)(v)\,s$. A general \(C\in|H|\) satisfies Petri's condition by \cite{lazarsfeld:brill-noether}. Hence \(\mu_T\) is injective and $\Ext^2(T, T) = 0$. Moreover, \(\operatorname{Hom}(T,T)=\mathbb{C}\), since \(L\) is a line bundle and \(V\hookrightarrow H^0(L)\). Consequently,
\[
\begin{aligned}
  \operatorname{ext}^1(T,T)&=1-\chi(T,T) = g+n(d+1-g-n) =\rho(g,n-1,d).
\end{aligned}
\]
The vanishing of \(\operatorname{Ext}^2(T,T)\) shows that the moduli space is smooth of this dimension at every point.

Therefore it is empty if \(\rho<0\). Classical Brill--Noether existence gives nonemptiness when \(\rho\ge 0\); for \(\rho=0\), smoothness and projectivity make it finite and reduced. For \(\rho>0\), Fulton--Lazarsfeld connectedness \cite{FL81}, together with smoothness, gives irreducibility.
\end{proof}

The main Theorem of this section is the following. 

\begin{Thm}\label{thm-rank one}
Let \(X\) be a smooth K3 surface with \(\operatorname{Pic}(X)=\mathbb{Z}H\), and let \(C\in|H|\) be a general smooth curve of genus \(g\ge 2\). Fix $d \in [1, 2g-2]$ and $n>1$. Then a general object $[\cO \otimes V \xrightarrow{\varphi} L] \in M_{\infty}(1, d, n)$ has surjective $\varphi$.
\begin{enumerate}
    \item If $n > d+1-g$, then every generated coherent system $[\cO \otimes V \xrightarrow{\varphi} L] \in M_{\infty}(1, d, n)$ has stable kernel $\ker \varphi$.
    \item If $n \leq d+1-g$, then a general generated coherent system $[\cO \otimes V \xrightarrow{\varphi} L] \in M_{\infty}(1, d, n)$ has stable kernel $\ker \varphi$.
\end{enumerate}
\end{Thm}

From now on, we fix $C \in |H|$ as in Theorem \ref{thm-rank one}. The following Lemma proves general generation.  

\begin{Lem} \label{lem-general}
    A general coherent system $[\cO \otimes V \xrightarrow{\varphi} L] \in M_{\infty}(1, d, n)$ has surjective $\varphi$ if $d \in [1, 2g-2]$ and $n>1$. 
\end{Lem}
\begin{proof}
Suppose that $T=[\mathcal{O}_C\otimes V\xrightarrow{\ \varphi\ }L]$ is not globally generated. Its image is $L(-D)$ for a unique nonzero effective divisor $D$, namely the common zero divisor of the sections in $V$. Thus $T$ fits into an exact sequence of coherent systems
\[
  0\longrightarrow\bigl[\mathcal{O}_C\otimes V\twoheadrightarrow L(-D)\bigr]\longrightarrow T\longrightarrow\bigl[0\longrightarrow L|_D\bigr]\longrightarrow 0.
\]
Conversely, a globally generated system $[\mathcal{O}_C\otimes V\twoheadrightarrow L']$ and an effective divisor $D$ determine $T$ up to isomorphism, by composing the evaluation map with the natural inclusion $L'\hookrightarrow L'(D)$.

Fix $\ell=\deg D>0$. By Lemma \ref{lem:K3-BN}, the globally generated systems of degree $d-\ell$ form a family of dimension at most $\rho(g,n-1,d-\ell)$, while $D$ varies in $\operatorname{Sym}^{\ell}C$, which has dimension $\ell$. Hence the locus of systems with base divisor of degree $\ell$, if nonempty, has dimension at most
\[
  \rho(g,n-1,d-\ell)+\ell=\rho(g,n-1,d)-(n-1)\ell<\rho(g,n-1,d).
\]
Since $1\le\ell\le d$, only finitely many such loci occur. By the same Brill--Noether lemma, each has dimension strictly smaller than every component of $M_\infty(1,d,n)$. A general system is therefore globally generated.
\end{proof}

Before proving part~(a) of Theorem~\ref{thm-rank one}, we recall a well-known lemma that transfers Brill--Noether information about line bundles on $C$ to higher-rank vector bundles.

\begin{Lem}\label{lem:gen-slope}
If $F$ is a globally generated vector bundle on $C$ of positive rank satisfying $H^0(F^\vee)=0$, then
\[
  \frac{\deg(F)}{\rk(F)}\ge 1+\frac{g}{\rk(F)+1}.
\]
\end{Lem}

\begin{proof}
Choose $\rk(F)+1$ general sections generating $F$. Dualizing their evaluation sequence gives
\[
  0\longrightarrow F^\vee\longrightarrow\mathcal{O}_C^{\,\rk(F)+1}\longrightarrow\det F\longrightarrow 0.
\]
Since $H^0(F^\vee)=0$, we have
\[
  h^0(\det F)\ge\rk(F)+1.
\]
Then Lemma \ref{lem:K3-BN} implies
\[
  g-\bigl(\rk(F)+1\bigr)\bigl(g-\deg(F)+\rk(F)\bigr)\ge 0.
\]
Rearranging gives $\deg(F)\ge\rk(F)+g\,\rk(F)/(\rk(F)+1)$, which is the assertion.
\end{proof}

\begin{Lem}\label{Lem:Rank one-above RR}
Fix $d \geq 1$ and $n>1$ such that $n > d+1-g$. Consider an arbitrary triple $[\cO \otimes V \xrightarrow{\mathrm{\varphi}} L] \in M_{\infty}(1, d, n)$ where $\varphi$ is 
surjective. Then $M:= \ker(\varphi)$ is slope-stable.
\end{Lem}

\begin{proof}
The evaluation sequence
\[
  0\longrightarrow M\longrightarrow \mathcal{O}_C\otimes V
  \longrightarrow L\longrightarrow 0
\]
gives \(\operatorname{rk}M=n-1\) and \(\deg M=-d\). If \(n=2\), then \(M\) is a line bundle, so the assertion is immediate.

Let \(F\subset M\) be a nonzero proper saturated subbundle, of rank \(a\le n-2\) and degree \(-e\). Since \(F\) is a subbundle of \(\mathcal{O}_C\otimes V\), its dual is globally generated. Moreover, $ H^0(F)=0$ because \(H^0(M)=0\), by the injectivity of \(V\to H^0(L)\). Then Lemma \ref{lem:gen-slope} for $F^{\vee}$ implies
\[
  \frac{e}{a}\ge 1+\frac{g}{a+1}.
\]
Thus using \(a+1\le n-1\) and \(d < n+g-1\), we obtain
\[
  -\mu(F) = \frac{e}{a}\ge 1+\frac{g}{a+1}\ge 1+\frac{g}{n-1} > \frac{d}{n-1} = -\mu(M),
\]
and so $M$ is slope-stable.
\end{proof}

\begin{proof}[Proof of Theorem \ref{thm-rank one}]
The general generation claim and part~(a) follow from Lemmas~\ref{lem-general} and~\ref{Lem:Rank one-above RR}, so it remains to prove part~(b). Recall that it assumes $n\leq d+1-g$. 

The case \(n=2\) is immediate since the kernel has rank one, we therefore assume $n \geq 3$. By Riemann--Roch, $h^0(L)\ge d+1-g\ge n$ for every $L\in\operatorname{Pic}^d(C)$. Hence the forgetful map
\[
  M_\infty(1,d,n)\longrightarrow\operatorname{Pic}^d(C),\qquad
  \bigl[\mathcal{O}_C\otimes V\to L\bigr]\longmapsto L,
\]
is surjective. By Serre duality, $L$ is special, i.e. $H^1(L)\ne 0$, precisely when $\omega_C\otimes L^{-1}$ is effective. The special locus is parametrized by effective divisors, thus has dimension at most $2g-2-d<g$, so a general line bundle of degree $d$ is nonspecial. Moreover,
\[
\rho(g,n-1,d)=g+n(d+1-g-n)\ge g>0,
\]
hence $M_\infty(1,d,n)$ is irreducible by Lemma \ref{lem:K3-BN}. It follows that a general coherent system in $M_\infty(1,d,n)$ has $H^1(L)=0$. 

Now choose the unique integer 
\[
  d_0\equiv d\pmod{n-1}\qquad \text{such that} \qquad  g+n-1\le d_0<g+2(n-1).
\]
By Proposition~\ref{prop:degree-shift-gen}, it suffices to find a single generated nonspecial system of degree $d_0$ with stable kernel. Indeed, repeated application of that proposition gives such a system in degree $d$, and openness of stability, together with the irreducibility of $M_\infty(1,d,n)$, proves the claim for a general system via the same argument as in Lemma \ref{lem-signle}.

Because of the choice of $d_0$ and the assumption $d\leq 2g-2$, we get
\begin{equation}\label{eq:d0-bounds}
  d_0\le g+2n-3,\qquad d_0\le 2g-2,\qquad n <g.
\end{equation}

We consider simultaneously the first destabilizing walls above the BN wall and proper equal-slope subbundles of the kernel at the BN wall. Take $T=[\cO \otimes V \overset{\varphi}{\twoheadrightarrow} L]\in \T$ with $\cl(T)=(1,d_0,n)$ and $H^1(L)=0$ which destabilizes at the wall $\mathcal{W}$. For $\mathcal{W}$ above the BN wall, retain the notation of Section~\ref{sec:wall-crossing}, namely
\[
T_1=[W\otimes\mathcal O_C\longrightarrow F],
\qquad \operatorname{cl}(T_1)=(a,e,k),
\qquad \mathcal{W}\colon \,w=tb+c,\qquad -\frac{1}{n-1}<b_0<0, 
\]
where $T_1$  is the first term of the Harder--Narasimhan filtration of $T$ just below $\mathcal{W}$. 

If a system $T$ reaches the BN wall, its kernel \(M_T\) is slope-semistable by Lemma \ref{lem-bn}. If $M_T$ is not slope-stable, choose a nonzero proper saturated subbundle
\[
M\subset M_T,\qquad
\mu(M)=-\frac{d_0}{n-1},
\]
and define
\[
W=\operatorname{im}\!\left(V^\vee\longrightarrow H^0(M^\vee)\right)^\vee
\subset V,
\qquad
F=(\cO \otimes W)/M, 
\]
thus $W$ is the smallest subspace with $M\subset  \cO \otimes W$, which implies that $F$ is a generated vector bundle with $H^0(F^\vee)=0$. This defines a subobject \(T_1 := [\cO \otimes W \to F]\hookrightarrow T\) in $\mathcal{A}(b_0)$ via a similar argument as in Step 3 of the proof of Theorem \ref{thm-any rank-high slope}.


Write $\cl(T_1)=(a,e,k)$, then from $\mu(M)=\mu(M_T)$ it follows
\[
\frac e{k-a}=\frac{d_0}{n-1},
\]
hence \(T_1\) has the same phase as \(T\) on the BN wall $\ell_f$ and is therefore semistable there.

Consequently, in both cases we have
\begin{equation}\label{eq:first-wall}
\begin{gathered}
k\le n,\qquad ad_0-e > 0,\qquad
n=c+td_0,\qquad k=ca+te,\\
c\ge1,\qquad
0\leq t\le\frac{n-1}{d_0}<\frac12.
\end{gathered}
\end{equation}
Indeed, since \(k-an=t(e-ad_0)\), the case of $t<0$ is impossible  by the first two inequalities in~\eqref{eq:first-wall}. Moreover, if $t=0$, the horizontal wall has \(a=1,\ k=n\). 


 The first claim is that \(a=1\) when $t>0$, i.e. all walls are made by objects from $\T$. Suppose otherwise, namely $a \geq 2$. By \cite[Lemma~6.5]{FN}, the vector bundle \(F\) is generated by $W$ with $H^0(F^\vee)=0$. Put \(M=\ker(\cO \otimes W\to F)\) when $\mathcal{W}$ lies above BN wall.

\medskip

\medskip
\noindent\emph{Step 1. Bounds on the rank and degree of $F$.} 
From~\eqref{eq:first-wall} we get 
$$
t(d_0-e)=n-k+c(a-1)\geq 1,
$$
therefore, together with an upper-bound for $t$ from~\eqref{eq:first-wall} and $d_0\leq 2g-2$, it gives the following upper-bounds on the degree of $F$:
\begin{equation}\label{f}
    e\leq \frac{n-2}{n-1}d_0\leq 2g-3.
\end{equation}

Next we bound the rank of $F$. The bundle \(M^\vee\) is generated with \(H^0(M)=0\), so applying Lemma \ref{lem:gen-slope} to $M^{\vee}$ gives
\[
  1+\frac{g}{k-a+1}\le\frac{e}{k-a} \overset{(*)}{\leq}\frac{d_0}{n-1} \overset{\eqref{eq:d0-bounds}}{\leq} 2+\frac{g-1}{n-1},
\]
where the inequality $(*)$ follows from~\eqref{eq:first-wall}, as
\[
  d_0(k-a)-(n-1)e=(c-1)(ad_0-e) \geq 0.
\]
Consequently, 
comparing the first term and the last term and using $n\geq 3$ we get
\begin{equation}\label{eq:n-vs-a}
  k-a+1 \geq \frac{g(n-1)}{g+n-2} \overset{\eqref{eq:d0-bounds}}{ > } \frac{n}{2},
  \qquad\text{hence}\qquad n\ge 2a-1,
\end{equation}
where the last one comes from $k\leq n$. This gives an upper-bound on the rank of $F$.

\medskip

\textit{Step 2. Excluding the case $a\geq 2$.} We do this by comparing the lower bound on the number of sections imposed by wall semistability with the upper bound coming from the Clifford bound for curves on K3 surfaces obtained in \cite{feyz-li:clifford-indices}. For this, let $F'$ be the last HN quotient of \(F\), write \(\rk(F')=r'\). It is slope-semistable, generated with $H^0(F'^\vee)=0$. 

We first claim that we can apply the next Lemma \ref{lem:clifford} to $F'$. We know
\begin{equation}\label{eq:mu-bound}
  0<\mu(F')\le\frac{e}{a} \overset{\eqref{f}}{\leq }\frac{2g-3}{a} < g-1.
\end{equation}
Let $W':= \text{im}(W \to H^0(F'))$ and $T_1'=[\cO \otimes W' \xrightarrow{\text{ev}|_{W'}} F']$. Then we have a surjection $T_1 \twoheadrightarrow T_1'$ in $\cA(b_0)$, where $(b_0, w_0)$ lies along the wall $\mathcal{W}$ made by $T_1$. Then 
$$
\nu_{b_0, w_0}([\cO \otimes W' \xrightarrow{\text{ev}|_{W'}} F']) = \frac{\frac{\dim W'}{\rk(F')} - w_0 }{\mu(F') -b_0} \geq \nu_{b_0, w_0}(T_1) = t
$$
which implies
\begin{equation*}
    \frac{h^0(F')}{\rk(F')} \geq \frac{\dim W'}{\rk F'} \geq w_0 + t(\mu(F')-b_0) = c+ t\mu(F').
\end{equation*}
In other words, $\Pi([\cO \otimes W' \xrightarrow{\text{ev}|_{W'}} F'])$ lies above the wall $\mathcal{W}$ and so $\ell_f$ with $\mu(F') \leq  \mu(F)$, thus 
\begin{equation}\label{eq:h0R-lower}
  \frac{h^0(F')}{r'}\ge \frac{\dim W'}{r'} \geq 1+\frac{n-1}{d_0}\,\mu(F').
\end{equation}
Hence to prove $h^0(F') \geq 2 r'$ as required in Lemma \ref{lem:clifford}, it suffies to show 
\begin{equation}\label{req}
    \frac{n-1}{d_0} \deg(F') > r' -1.
\end{equation}
Applying Lemma \ref{lem:gen-slope} to $F'$ gives
\begin{equation*}
    \deg(F') \geq r' + \frac{g r'}{r'+1} \overset{(*)}{>} \frac{g}{2} + 2(r'-1) \overset{(**)}{ \geq} \frac{d_0(r'-1)}{n-1},
\end{equation*}
where $(*)$ comes from 
\[
  r'+\frac{gr'}{r'+1}-\left(\frac{g}{2}+2(r'-1)\right)
  =\frac{(r'-1)(g-2r')+4}{2(r'+1)}>0, 
\]
where the last strict inequality is the result of $g>n$ from \eqref{eq:d0-bounds}, $n \geq 2a-1$ from \eqref{eq:n-vs-a} and $a \geq r'$. On the other hand, $(**)$ is equivalent to 
$$
\frac{g}{2} \geq \frac{r'-1}{n-1} \left( d_0 -2 (n-1) \right), 
$$
which holds since $g> d_0 -2(n-1)$ by \eqref{eq:d0-bounds} and $\frac{n-1}{2} \geq a -1 \geq r' -1$ as above. Thus \eqref{req} holds, and consequently $h^0(F')\geq 2 r'$. 

In order to apply Lemma \ref{lem:clifford} to $F'$, it remains to show $(r')^2<g$. By~\eqref{eq:mu-bound}, we have $\mu(F')\le(2g-3)/a$, thus the quadratic section bound \cite[Theorem~1.1]{feyz-li:clifford-indices} gives
\[
  2\le\frac{h^0(F')}{r'}<1+\frac1g+\frac{g(2g-3)^2}{4a^2(g-1)^2}.
\]
Rearranging and using $r'\le a$, we obtain
\[
  (r')^2\le a^2<\frac{g^2(2g-3)^2}{4(g-1)^3}<g,
\]
where the last inequality follows from $4(g-1)^3-g(2g-3)^2=3g-4>0$.

\medskip

Thus Lemma \ref{lem:clifford} applies to $F'$. Since $\max\{r', 2\}\leq a$, it gives 
\begin{equation*}
  \frac{h^0(F')}{r'}\le 1+\frac{\mu(F')}{2}-\frac{g-2}{2a}.
\end{equation*}
Comparing with \eqref{eq:h0R-lower}, we obtain 
$$
\frac{\mu(F')}{2} - \frac{g-2}{2a} \geq \frac{n-1}{d_0}\mu(F') \overset{\eqref{eq:d0-bounds}}{\geq} \frac{\mu(F')}{2+ \frac{g-1}{n-1}}. 
$$
Rearranging and applying Step~1 alongside with $n<g$ and $d_0\leq g+2n-3$ from~\eqref{eq:d0-bounds}, we obtain
\[
  \frac{(g-2)(g+2n-3)}{g-1}\le a\,\mu(F')\le e \overset{\eqref{f}}{\leq }\frac{n-2}{n-1}\,d_0
  <\frac{(g-2)(g+2n-3)}{g-1}.
\]
This contradiction proves that $a=1$.
\medskip

\textit{Step 3. Dimension count.} It remains to count the dimension loss along the wall $\mathcal{W}$. First, $a=1$ implies that $f\colon F \to L$ is injective, therefore we can write
\[
F=L(-D),\qquad \ell=\deg D=d_0-e>0,
\]
for an effective divisor $D$. In particular, the destabilizing sequence for $T$ is the following exact sequence in $\T$
\[
  0\longrightarrow\bigl[\mathcal{O}_C\otimes W\twoheadrightarrow L(-D)\bigr]\longrightarrow T\longrightarrow\bigl[\cO\otimes (V/W) \to  L|_D\bigr]\longrightarrow 0.
\]
The case $k=n$ is excluded by the surjectivity of $\varphi$. Moreover, $k\geq 2$ since $F$ is generated and $H^0(F^\vee)=0$. Thus we can assume $2\le k\le n-1$. 

Conversely, every such $T$ arises from a choice of a generated system $[\cO \otimes W \twoheadrightarrow L']$ of class $(1, e, k)$ together with an effective divisor $D$ of degree $l$ and a choice of an $n$-dimensional subspace $V$ of $H^0(L)$ containing $W$. By Lemma~\ref{lem:K3-BN}, the system $[\cO \otimes W \twoheadrightarrow L']$ varies in a family of dimension at most $\rho(g, k-1, e)$, while $D$ varies in $\text{Sym}^l(C)$, which has dimension $l$. Since $L$ is nonspecial, we have that $h^0(L)=d_0+1-g$, thus the remaining choices form $\text{Gr}(n-k, H^0(L)/W)$. 

Hence the locus of nonspecial systems that gets destabilized along $\mathcal{W}$ has dimension at most
\[
\rho(g,k-1,e)+\ell+(n-k)(d_0+1-g-n).
\]
Using $l=d_0 -e$, its codimension in $M_{\infty}(1, d_0, n)$ is consequently at least
\[
(k-1)(d_0-e) -k(n-k) \overset{(*)}{\geq} (n-k)e -k(n-k)=(n-k)(e-k), 
\]
where $(*)$ follows from~\eqref{eq:first-wall} as $k-1=(c-1)+te\geq te$ and $n-k=t(d_0-e)$.  
Furthermore, Lemma \ref{lem:K3-BN} gives $\rho(g, k-1, e)\geq 0$, equivalently 
\[
e\ge(k-1)\left(1+\frac gk\right)>2(k-1)\ge k,
\]
since \(g>n>k\) and $k\geq 2$. Thus, the codimension is strictly positive. Only finitely many pairs \((e,k)\) occur. A general system therefore avoids every earlier wall and every equal-slope locus on the final wall, so its kernel is slope-stable. This proves part (b).

\end{proof}

\begin{Lem}\label{lem:clifford}
Let \(E\) be a slope-semistable vector bundle on \(C\) of rank \(r\) with \(r^2\le g\), satisfying
\[
  \mu(E)\le g-1,\qquad h^0(E)\ge 2r,
\]
then
\[
  \frac{h^0(E)}{r}\le 1+\frac{\mu(E)}{2}-\frac{g-2}{2\max\{r,2\}}.
\]
\end{Lem}

\begin{proof}
If \(r\ge 2\), the Clifford bound of \cite{feyz-li:clifford-indices} applies because \(g\ge r^2\), giving
\[
  \mu(E)-\frac{2h^0(E)}{r}+2
  \ge\frac{2}{r}\left(g-1-\left\lfloor\frac{g}{r}\right\rfloor\right)
  \ge\frac{g-2}{r}.
\]
If \(r=1\), then \(h^1(E)\ge h^0(E)\ge 2\). The ordinary Clifford bound gives
\[
  \mu(E)-2h^0(E)+2\ge\operatorname{Cliff}(C)=\left\lfloor\frac{g-1}{2}\right\rfloor\ge\frac{g-2}{2}.
\]
Rearranging proves the claim.

\end{proof}

We finish this section with the proof of Corollary \ref{cor:bundles} in the Introduction. 

\begin{proof}[Proof of Corollary \ref{cor:bundles}]

Suppose first that \(F\) is a globally generated semistable bundle of rank \(R\) and degree \(d>0\). Since \(F^\vee\) is semistable of negative slope, \(H^0(C,F^\vee)=0\). As in the proof of Lemma~\ref{lem:gen-slope}, choose \(R+1\) general sections generating \(F\). Such a choice is possible because failure to generate at a fixed point has codimension two in \(H^0(C,F)^{\oplus(R+1)}\), while the point varies on a curve. Dualizing the resulting evaluation sequence gives
\[
0\longrightarrow F^\vee\longrightarrow\mathcal O_C^{\oplus(R+1)}\longrightarrow\det F\longrightarrow0.
\]
It follows that \(h^0(C,\det F)\ge R+1\). If \(d\le2g-2\), then Lemma \ref{lem:K3-BN} gives
\[
\rho(g,R,d)=g-(R+1)(g-d+R)\ge0.
\]
If \(d>2g-2\), Riemann--Roch instead gives
\[
R+1\le h^0(C,\det F)=d+1-g,
\]
and hence \(R\le d-g\). This proves necessity.

Conversely, suppose that the stated numerical condition holds. If \(d\le2g-2\), Lemma \ref{lem:K3-BN} gives a line bundle \(L\) of degree \(d\) with at least \(R+1\) sections. If \(d>2g-2\), every line bundle \(L\) of degree \(d\) has \(h^0(C,L)=d+1-g\ge R+1\). Thus \(M_\infty(1,d,R+1)\) is nonempty in either case. By Theorems~\ref{thm-rank one} and~\ref{thm-any rank-high slope}, a general system in this space has an exact evaluation sequence
\[
0\longrightarrow M_{L,V}\longrightarrow V\otimes\mathcal O_C\longrightarrow L\longrightarrow0
\]
with \(M_{L,V}\) slope-stable. Dualizing gives
\[
0\longrightarrow L^\vee\longrightarrow V^\vee\otimes\mathcal O_C\longrightarrow M_{L,V}^\vee\longrightarrow0.
\]
Hence \(M_{L,V}^\vee\) is globally generated, has rank \(R\) and degree \(d\), and is slope-stable. In particular, the numerical conditions suffice even for the existence of a stable bundle.


\end{proof}

\bibliography{mybib}
 \bibliographystyle{halpha}
 
\end{document}